\documentclass[12pt]{amsart}
\usepackage{amsmath,amsthm,amsfonts,amssymb,verbatim}
\usepackage{hyperref}
\usepackage[usenames]{color} 
\usepackage{graphicx}
\usepackage{tikz,tikz-cd}
\usetikzlibrary{arrows}
\usetikzlibrary{arrows.meta,decorations.pathreplacing,patterns}

\usepackage{mnotes}
    
\definecolor{teal}{RGB}{22, 180, 180}

\newtheorem{thm}{Theorem}[section]
\newtheorem{theorem}[thm]{Theorem}
\newtheorem{lemma}[thm]{Lemma}
\newtheorem{cor}[thm]{Corollary}
\newtheorem{proposition}[thm]{Proposition}
\newtheorem{prop}[thm]{Proposition}

\newtheorem{question}[thm]{Question}

\usepackage{pdfsync}
\usepackage{graphicx}

\theoremstyle{definition}  

\numberwithin{equation}{section}
\newtheorem{defn}[thm]{Definition}

\newtheorem{ex}[thm]{Example}

\newtheorem{example}[thm]{Example}
\newtheorem{obs}[thm]{Observation}
\newtheorem{remark}[thm]{Remark}

\theoremstyle{definition}
\theoremstyle{remark}

\newcounter{enumitemp}

\DeclareMathOperator{\card}{card}

\DeclareMathOperator{\Int}{int}

\DeclareMathOperator{\lo}{o}

\DeclareMathOperator{\diam}{diam}

\newcommand{\dsp}{{\delta^+}}

\newcommand{\fB}{{\mathfrak B}}

\newcommand{\B}{{\mathcal B}}
\newcommand{\C}{{\mathcal C}}

\newcommand{\eps}{\varepsilon}
\newcommand{\cA}{{\mathcal A}}
\newcommand{\cB}{{\mathcal B}}
\newcommand{\cD}{{\mathcal D}}
\newcommand{\E}{{\mathcal E}}
\newcommand{\F}{{\mathcal F}}
\newcommand{\cH}{{\mathcal H}}
\newcommand{\cI}{{\mathcal I}}
\newcommand{\cK}{{\mathcal K}}

\newcommand{\U}{{\mathcal U}}
\newcommand{\V}{{\mathcal V}}
\newcommand{\R}{\mathbb R}
\newcommand{\fF}{{\mathfrak F}}

\newcommand{\N}{{\mathbb N}}

\newcommand{\cP}{{\mathcal P}}

\newcommand{\cS}{{\mathcal S}}

\newcommand{\Z}{\mathbb Z}

\def\cG{{\mathcal G}}

\def\cS{{\mathcal S}}

\def\Q{{\mathbb Q}}

\def\cT{{\mathcal T}}
\def\z2s{{$\Z^2$-subshift}}

\def\CA{{\mathcal A}}

\def\Z{{\mathbb Z}}

\DeclareMathOperator{\End}{{\rm End}}

\title{Polygonal $\Z^2$-subshifts}
\author{John Franks}
\address{Northwestern University, Evanston, IL 60208 USA}
\email{j-franks@northwestern.edu}
\author{Bryna Kra}
\address{Northwestern University, Evanston, IL 60208 USA}
\email{kra@math.northwestern.edu}

\subjclass[2010]{}
\keywords{subshift, automorphism, nonexpansive}

\thanks{The second author was partially supported by NSF grant DMS-1500670.}

\begin{document}

 \begin{abstract}
 Let  $\cP\subset\Z^2$ be a convex polygon with each vertex in it labeled by an element from a finite set 
 and such that the labeling of each vertex $v\in \cP$ is uniquely determined by the labeling of all 
 other points in the polygon.  We introduce a class of $\Z^2$-shift systems, the {\em polygonal shifts}, determined by such a polygon: these are shift systems such that the restriction of any $x\in X$ to some polygon $\cP$ has this property. 
 These polygonal systems are related to various well studied classes of shift systems,  
 including subshifts of finite type and algebraic shifts, but include many other systems.  
We give necessary conditions for a $\Z^2$-system  $X$ to be polygonal, in terms of the nonexpansive subspaces of $X$, 
and under further conditions can give a complete characterization for such systems.  
  \end{abstract}
\maketitle

\section{Introduction}
If $\cA$ is a finite alphabet, a {\em $\Z^2$-shift} $X$ is a closed subspace of $\cA^{\Z^2}$ that is invariant 
under the $\Z^2$-action by horizontal and vertical shifts.  
Large classes of shifts have been well studied, including algebraic shifts and
shifts of finite type (see for example~\cite{schmidt, kitchens, LM}).  We focus on a collection related to
these, which we call polygonal shifts.

Roughly speaking, polygonal shifts are a class in which the data in one
region determines the data in a larger region.  
We defer the precise definitions
until Sections~\ref{sec:basic-defs} and~\ref{sec:polygonal}, starting with some examples that motivate the study of these shifts.  
We refer to an element $x = \bigl(x(i,j)\colon i,j\in\Z\bigr)$ in a $\Z^2$-shift $X$ as {\em coloring}, and 
refer to the restriction of $x\in X$ to a set $S\subset\Z^2$ as a {\em coloring of $S$}.  
Perhaps the
simplest interesting example is the 
Ledrappier shift~\cite{ledrappier}: if $\CA = \Z/2\Z$, 
define $X$ to be the subshift of $\CA^{\Z^2}$ such that every $x\in X$ satisfies 
\begin{equation}
\label{eq:ledrappier}
x(i,j) + x(i+1, j) + x(i, j+1) = 0 \mod 2.  
\end{equation}
The key property is that for the triangle $\cT$ with vertices 
$(0,0), (0,1)$, and $(1,0)$, 
the coloring of a vertex is uniquely determined
by the coloring of the other two vertices of the triangle, and this triangle is what motivates 
the commonly used name {\em three dot system} for this shift.   Shift invariance of the system
implies that the same holds for vertices of any integer translate of $\cT$.  

Polygonal systems generalize this idea, and instead of using a triangle, 
we consider an arbitrary convex polygon $\cP\subset \R^2$, 
and we refer to such a polygon as an {\em (integer)
  polygon} if all its vertices lie in $\Z^2$.
The key property of the polygon 
$\cP$ is that the color of each vertex $v\in\cP$ is uniquely determined by the coloring of all points of $\cP$ 
other than $v$ (note that these other points may include interior points of the polygon).
If $X$ is a $\Z^2$-shift and there is a convex polygon
with vertices in $\Z^2$ such that the restrictions of all colorings $x\in X$ 
to the polygon $\cP$ has this property, then we say that the system is
{\em polygonal with respect to $\cP$} and that $\cP$ is a {\em coding
  polygon} for the shift.  We emphasize that by definition coding polygons
have their vertices in $\Z^2$ and hence all edges have rational slopes.

One of our goals is to characterize the $\Z^2$-shifts which are
polygonal and for a polygonal shift ascertain to what extent its
coding polygon $\cP$ is canonical.  In general, coding polygons are
far from unique.  If a shift is polygonal with respect to
$\cP$, then it is obviously also polygonal with respect to
$\cP +(i,j)$ for any $i,j \in \Z$.  Furthermore, it is also polygonal 
for the convex hull of
$\cP$ and $\cP +(i,j)$, and this hull is a coding polygon with two
additional sides parallel to $(i,j)$ (assuming there is not already a
side of $\cP$ parallel to this vector).  Repeating this construction
shows that coding polygons for a given shift can have an arbitrarily
large number of sides and that any finite set of rational slopes can
be among the slopes of their edges.  In light of this it is natural to
ask what is the minimal number of edges of a coding polygon and what
edge slopes must occur in any coding polygon for $X$. We address these
questions in Theorem~\ref{main thm1} and Proposition~\ref{prop:
  polygon sides}, respectively. 

The key concepts in addressing these and other
questions are the notions of expansive and
nonexpansive.
 Again, we postpone the formal definitions until Section~\ref{sec:basic-defs}, but 
we motivate their role.  For the 
Ledrappier system $X$, it is easy to check that for all but three (up
to translation) half spaces in $\R^2$, any coloring of its integer
lattice points extends uniquely to a coloring of all of $\Z^2$; this is
well known, and follows from a more general result given in
Proposition~\ref{prop: polygon sides}.  The only exceptions are the
half spaces which are translates of the three half spaces
given by the inequalities  $x \ge 0$, \ $y \ge 0$, and 
$y \le -x$.
To make precise the sense in which data in one region determines data
outside this region, we view a half space as being specified by an
oriented ray, namely a ray which lies in the boundary of that half
space and inherits its orientation from the induced orientation on the
boundary.  If every coloring of the integer lattice in a
half space extends uniquely to the full
space, we say that the half space and corresponding oriented ray are
{\em expansive}, and otherwise we say that they are {\em nonexpansive}.  
In the Ledrappier system, the only nonexpansive rays are  the rays
lying in the boundaries of the three specified half spaces and
having the appropriate orientation, namely the rays
spanned by the vectors $(1,0), (-1,1)$ and $(0, -1)$.
Note that these three vectors form the edges (not vertices) of an
oriented coding polygon for the Ledrappier shift.

This terminology is consistent with standard notions of expansiveness
and nonexpansiveness for one-dimensional subspaces of $\Z^2$-systems, as 
studied, for example, in Boyle and Lind~\cite{BL}.  In particular, they show that
the set of nonexpansive subspaces is nonempty when $X$ is infinite.
Accordingly, any $\Z^2$-shift $X$ we consider is assumed to
be infinite.
Allowable colorings of a nonexpansive half
space do not uniquely determine the coloring of even a single point in
the complementary half space, and this behavior again shows up in the Ledrappier system.  
In our more general setting of polygonal systems, it is exactly the nonexpansive rays that are 
used to characterize which shifts lie in this class. 

Generalizing the Ledrappier example,  Kitchens and Schmidt~\cite{KS, KS2}
study $\Z^d$-actions on Markov subgroups.  If $\cA$ is a  finite
abelian group $\cA$, then $\cA^{\Z^d}$ is a zero-dimensional compact abelian
group when endowed with the operation of component-wise addition. A {\em Markov subgroup} $X$ is a closed
subgroup of this group such that there exists some finite set
$S \subset \Z^d$ (called a {\em shape}) satisfying 
\[
\sum_{u \in (S+v)} x(u) = 0
\]
for each fixed $v \in \Z^d$.  
For $d =2$, it is easy to see that any Markov subgroup is polygonal with
coding polygon $\cP$ given by the convex hull of the finite set $S$.

The polygonal shifts are a class of zero-dimensional
$\Z^2$-subshifts that is more general and substantially larger than Markov subgroups or
similar systems with a strong algebraic structure.  More precisely, a
result of Einsiedler~\cite{E} shows there are uncountably many
$\Z^2$-invariant subspaces of the Ledrappier shift $X$ with distinct
topological entropies and it follows that there are uncountably many distinct polygonal
shifts with the same alphabet and the same polygon $\cT$. In
particular, not all polygonal shifts are isomorphic to
subshifts of finite type or to $\Z^2$-actions on Markov
subgroups, as these classes are countable  (up to isomorphism).

We  limit ourselves to shifts which are polygonal with respect to convex polygons.  There
is no loss in doing so, as it is easy to see that a shift which is polygonal with respect to
a polygon $\cP$ is also polygonal with respect to the convex hull $\hat \cP$ of $\cP$. 
The advantage of working with $\hat \cP$ is that it has strictly fewer edges and vertices than $\cP$, 
unless $\cP$ is already 
convex.

Another reason to make use of the simplification in the geometry in
passing to the convex hull of a shape, rather than more general
shapes, is that the edges of a convex coding polygon are closely
related to the geometry of nonexpansive subspaces.  For example, if
$X$ is a Markov subgroup with shape $S$, then the nonexpansive
subspaces are precisely the subspaces parallel to the edges of $\cP$, the
convex hull of $S$.

The fact that all other subspaces are expansive is a special case of a result given in 
Proposition~\ref{prop: polygon sides}.  
Since for each edge $e$
there are multiple legal colorings of $\hat \cP$ which differ on $e$
but agree on $\hat \cP \setminus e$, it follows that the 
 edges are 
nonexpansive (see Definition~\ref{def expansive}).

In seeking the simplest polygon to represent a shift $X$ we  allow
ourselves to replace $X$ with a particular kind of isomorphic shift $Y$ which we call a {\em
  recoding} of $X$.  The precise definition of recoding is given in
Definition~\ref{def: recoding}, but again we give an informal motivation.  
Starting with a finite convex subset $F\subset\Z^2$, we create a new alphabet $\cA_F$ consisting 
of all legal colorings of $F$.  The recoding
$X_F$ of $X$ is then the subset of $(\cA_F)^{\Z^2}$ with the property that for each $y \in X_F$
there is an $x \in X$ such that for each $i,j \in \Z^2$, the coloring of $y(i,j)$ is the restriction
of the coloring $x$ to $F + (i,j)$.

Considering the class $\cP(X)$ of
all (integer) coding polygons for all recodings of a subshift $X$, we refer to a
polygon $\cP_0 \in \cP(X)$ as a {\em minimal recoding polygon} if it has
the minimal number of edges of all elements of $\cP(X)$ and
 is minimal under inclusion among coding polygons with
  that number of edges.
Note that a minimal recoding polygon for $X$ is
a coding polygon for a recoding of $X$, not necessarily for $X$ itself.

We show in Proposition~\ref{scale down} that if
a coding polygon $\cP$ for $X$ is equal to $n \cP_0$ for some integral polygon
$\cP_0$, then $X$ can be recoded to a polygonal system with a coding
polygon $\cP_0$.  Hence a minimal recoding polygon must be primitive in
the sense that it is not an integer multiple of a smaller integer polygon.
A natural question arises: what are the possible minimal recoding polygons
for a polygonal shift? 

The geometry of  the minimal recoding polygons for
a polygonal system $X$ is closely linked
to the nonexpansive rays of $X$.
To make this more precise we refine the notion
of parallel to distinguish whether parallel objects have orientations
which coincide (see Section~\ref{sec:parallelisms} for complete definitions).   
We view a {\em ray} in $\R^2$ as
a translate of the set $\ell_v = \{tv \colon t \ge 0, \ v\ne 0 \in \R^2\}$, 
and assume it is oriented in the direction of increasing $t$.  
We refer to two rays $\ell_1$ and $\ell_2$  which are translates of
  each other  as  {\em positively parallel}, and when 
the rays  $\ell_1$ and $-\ell_2$ are positively parallel, 
we say $\ell_1$ and $\ell_2$  are {\em antiparallel}.  Thus the standard understanding of rays being 
parallel means they are either positively parallel or antiparallel.
We extend these conventions to oriented line segments, referring to such as a segment 
as positively parallel  to a ray $\ell$ (or to another line segment) if it has the same orientation and otherwise 
as antiparallel to $\ell$ (or again to another line segment).

Orientations extend to polygons $\cP\subset \R^2$: such a polygon 
inherits an orientation from $\R^2$, and this orientation induces an
orientation on the boundary $\partial \cP$ and hence an orientation on
each edge of $\cP$.  
In a convex polygon, no two edges can be positively parallel, but pairs of edges may
be antiparallel.

\begin{thm}
If $X$ is an infinite polygonal shift and $\cP_0$ is a minimal
recoding polygon for some recoding of $X$, then any ray positively parallel to
the oriented edges of $\cP_0$ is nonexpansive for $X$ and every
other ray is expansive.
\end{thm}

This result is an immediate consequence of Proposition~\ref{prop: polygon sides} and
Theorem~\ref{main thm1} and Theorem~\ref{main thm2}, and it provides a necessary
condition for a $\Z^2$-system to be polygonal; it must have finitely many
nonexpansive rays and they must all have rational slope.

However, an example
of Hochman~\cite{hochman} shows that this condition is not sufficient. 
There is an additional necessary property, called {\em closing}  (see 
Definition~\ref{def closing}), that must 
be satisfied by the nonexpansive rays in polygonal systems. With this
additional hypothesis we have both necessity and sufficiency:

\begin{thm}
Suppose $X$ is an infinite $\Z^2$-subshift with finitely many nonexpansive rays 
each of which has rational slope  and  is closing.
Then there is a recoding $Y$ of $X$ which is polygonal.
\end{thm}

This result follows from Theorem~\ref{main thm1}, and in
Theorem~\ref{main thm2}, we give a version of the converse: if $\cP_0$ is a coding polygon
for $X$, then $X$ can be recoded to a
subshift $Y$ with a coding polygon $\cP$ having $m$ edges, the minimum
possible.

Moreover any two such minimal recoding polygons have parallel edges
(and hence have equal corresponding angles).  If $\cP_0$ is a
triangle, we can say more and in Corollary~\ref{cor: min triangle}, we
show that if $\cP_0$ is a minimal recoding triangle for an infinite
$X$, then it is uniquely determined up to translation. We do not know
if this generalizes, and in particular do not know if minimal recoding polygons
which are not triangles are unique up to translation.

A system isomorphic to a polygonal system
need not be polygonal, but in Corollary~\ref{cor: recoding poly} we show
that if $Y$ is a recoding of $X$ and
$X$ is polygonal then  so is $Y$.

In Section~\ref{sec:entropy}, we study various forms of entropy for
$\Z^2$-systems.  For an arbitrary $\Z^2$-system $X$ and direction
$v\in \Z^2$, there is a seminorm $\|\cdot \|_X$ that captures the
directional entropy for $X$ in direction $v$ (see~\cite{BL} and~\cite{M}). 
In Corollary~\ref{cor: poly norm} we observe that a result of Milnor implies
that for any polygonal
system, whose coding polygon has no antiparallel 
sides, this seminorm $\|\cdot \|_X$ is either identically zero or is a
norm.  Furthermore, in Proposition~\ref{prop quasi-conformal} we show
that in this case, if the entropy norms are nontrivial, then
the associated seminorms for the family of polygonal systems
associated to a given polygon is a quasi-conformal family.
  Roughly speaking, this means that
for any subshift $Y$ in the same family as a subshift $X$, a sphere in
the norm $\|\cdot \|_X$ has bounded eccentricity in $\|\cdot \|_Y$
with a bound that is independent of $Y$.

More precisely, suppose
$\cP$ is a  rational polygon 
which has no antiparallel edges and 
$\fF(\cP)$ is the family of all $\Z^2$-subshifts which are polygonal with respect to
$\cP$ and which have nontrivial entropy norms, 
Then we show (Proposition~\ref{prop quasi-conformal}) that 
there is a uniform dilatation constant $D>0$, depending only on $\cP$, which
has the property that for all $X \in \fF(\cP)$ and
any $u, v \in S^1$  we have
\[
\frac{1}{D} \le \frac {h_u(X)}{h_v(X)} \le D.
\]

When the polygon is a triangle we obtain a stronger result, showing that 
they are conformally equivalent.  

In Corollary~\ref{cor conformal}, we show that 
if $X, Y$ are triangular $\Z^2$-subshifts  with  nontrivial entropy norms
with respect to the same rational triangle $\cT$, then
there is a constant $C >0$ such that $\|\cdot \|_X = C \|\cdot\|_Y $  and 
the constant does not depend on the direction  chosen in $\R^2$.

\subsection*{Acknowledgment}
We thank Van Cyr for many invaluable conversations during the preparation of this paper
and we thank the referee for numerous comments that improved the paper.

\section{Background on shift systems}
\label{sec:basic-defs}
\subsection{Shift systems}
We assume throughout that $\cA$ is a finite set, called the {\em alphabet}, endowed with the discrete
topology.  For $d\geq 1$, we endow $\cA^{\Z^d}$ with the product topology. We
review the standard definitions for $\cA^{\Z^d}$ for any $d\geq 1$ when there is no notational difference, but in most of the article we focus on two dimensions.

 An element $x\colon \Z^d \to \cA$ is called a {\em coloring}
and $x(u)$ denotes the color of $x$ at the position $u\in\Z^d$.
When we want to make use of both coordinates in two dimensions, we use $x(i,j)$ to denote
the color of $x$ at the position $(i,j)\in\Z^2$.

If $X\subset \cA^{\Z^d}$
is closed and invariant under the $\Z^d$ action $(T^{u}x)(v) = x(u+v)$ for $u\in\Z^d$,
then we say that $X$ is {\em $\Z^d$-subshift},
and when the context is clear, we shorten this and say that $X$ is a {\em shift system} or just a {\em shift}, omitting the transformations from the notation.  Thus
in two dimensions, such $X$ is implicitly endowed with the horizontal $T^{(1,0)}$
and vertical $T^{(0,1)}$ shifts.
When  considering more than one shift possibly with different alphabets,
we write $(X, \cA)$ to emphasize the alphabet and, by convention, we only include in
$\cA$ letters which are used in the language of $X$.
If in addition we need to distinguish the transformations on different shifts, we write $(X, T_X)$, or
$(X, \cA, T_X)$  when we need to capture all of the data.
When there is no possible ambiguity,
we refer to transformations $T_X$ and $T_Y$ on different spaces $X$ and $Y$ as just $T$.

If $X$ is a subshift, we refer to an element $x\in X$ as an {\em $X$-coloring} and we refer to the restriction of $x$ to a region $A$ as an $X$-coloring of $A$. 

\subsection{Coding and Recoding}
\label{sec:recoding}

Of particular interest is how the coloring information from one region
in $\Z^d$ forces the coloring of another region, or perhaps all of $\Z^d$.  We recall a definition from~\cite{BL} which makes this precise:
\begin{defn}
If  $X$ is a $\Z^d$-subshift and
$A,B \subset \R^d$, then $A$ {\em $X$-codes} $B$
if for all $x, x' \in X$, whenever $x$ and $x'$
agree on $A \cap \Z^d$, then they also agree on $B \cap \Z^d$.
If the shift $X$ is clear from the context, we just say that $A$ {\em codes} $B$.
In a slight abuse of notation, we say $A$ codes $v \in \Z^d$ to
mean that $A$ codes the set $\{v\}$ of a single element.
\end{defn}

Note that there is no assumption that the region $A$ is finite, and the definition is stated for $A, B$ as subsets of $\R^d$.
Though the configurations $x, x'\in X$ are only defined on integral coordinates, the more general definition of the subsets gives us necessary flexibility for some of the results.

A trivial example of a region coding another is in a doubly periodic shift, where any set $A$ that contains a full period completely determines an entire configuration and so codes all of $\cA^{\Z^2}$.  At the opposite extreme is the full shift $\cA^{\Z^2}$;
no region codes any larger region.

Since by definition a shift system is translation invariant, we have the following
immediate fact:

\begin{remark}\label{rem:translate}
Since a shift $X$ is invariant under the $\Z^d$-action, it follows immediately that
for every $v \in \Z^d$, if $A$ codes $B$ then $A + v$ codes
$B + v$.
\end{remark}

Recall that an {\em isomorphism} $\Psi\colon (X, \cA, T_X)\to (Y, \cA', T_Y)$ is homeomorphism
such that $\Psi\circ T_X = T_Y \circ \Psi$.

\begin{defn}\label{def: recoding}
If $(X, \cA)$ is a $\Z^2$-shift and $F$ is a finite subset
of $\Z^2$, we say the $\Z^2$-shift $(Y, \cA')$ is a {\em recoding of $(X, \cA)$ via $F$}
provided there is an
isomorphism of $\Z^2$-shifts $\Psi\colon (X, \cA) \to (Y, \cA')$
such that for every $(i,j) \in \Z^2$ and all $x,x' \in X$
\[
\Psi(x)(i,j) = \Psi(x')(i,j) \ \text{ if and only if  }\
x|_{F(i,j)} = x'|_{F(i,j)}
\]
where $F(i,j) := F + (i,j)$.
Equivalently $\{(i,j)\}$ \ $\Psi^{-1}$-codes $F(i,j)$ and
$F(i,j)$\ $\Psi$-codes $\{(i,j)\}$.

Note that we slightly overload notation but it should be clear from the context what is meant.
We use capital letters such as $F$ or $R$ for
subsets of $\Z^2$ and in this case,
for example, $F(i,j)$ denotes the translate of the set
$F + (i,j) = \{f + (i,j)  \colon f \in F\}$ , while we use lower case letters
such as $x$ or $y$ for elements of a shift $X$ and in this
case, for example,  $x(i,j)$ denotes the color which $x$  assigns to $(i,j)$.

Given  $X$ and any finite subset $F \subset \Z^2$,
we define the
{\em canonical recoding $(X_F, \cA_F,)$ via $F$} by setting
$\cA_F$ to be the set of all colorings of $F$ which are the restriction of
colorings in $X$ and setting $\Psi$ to be the isomorphism induced by the
block map which assigns to a restriction to $F$ of an
$X$-coloring the element of $\cA_F$ it represents.
\end{defn}

Recall that a map $\Psi\colon X\to Y$ is an {\em $r$-block code} if for all $x\in X$,
the color that $\Psi(x)$ assigns to $0$ is determined by the values of $x(i,j)$
with $\| (i,j) \| \le r$ (when needed, we use the Euclidean norm $\|\cdot\|$ on $\R^2$).

If $\Psi\colon (X, \cA) \to (Y, \cA')$ is a recoding, then its inverse is an isomorphism
induced  by a $0$-block map $\phi\colon  \cA' \to \cA$.
It is clear that the relation ``$Y$ is a recoding of $X$'' is reflexive.
It is also transitive, because the composition of two recodings is
a recoding.  However, this relation is not symmetric, as
whenever $(Y, \cA')$ is a recoding of $(X, \cA)$, it follows that
$\card(\cA') \ge \card(\cA)$ and this inequality is usually strict.
Indeed if $Y$ is a recoding of $X$ and $X$ is a recoding of $Y$,
then there is a bijection of their respective alphabets which induces
an isomorphism.

We now show if $(Y, \cA')$ is a recoding of $(X,\cA)$ via a
finite set $F$, then there is an isomorphism of $Y$ with $X_F$ induced
by a bijection of $\cA'$ and $\cA_F$.

\begin{lemma}\label{lem:recoding}
Suppose $X$ is a $\Z^2$-shift, $F \subset \Z^2$ is finite,
and $X_{F}$ is the canonical recoding.
\begin{enumerate}
\item If $v \in \Z^2$ and $T_X^v$ is the shift on $X$ corresponding to $v$,
then $\Psi\colon  (X, \cA) \to (Y, \cA')$ is a recoding via $F$ if and only if
$T_Y^v \circ \Psi$ is a recoding of $X$ via $T_X^v(F)$.

\item If $(Y, \cA')$ is a recoding of $(X, \cA)$ via ${F}$, then there is an isomorphism
of $(Y, \cA')$ and the canonical recoding $(X_{F}, \cA_{F})$ induced by a bijection of
the alphabet $\cA'$ with the alphabet $\cA_{F}$.
\end{enumerate}
\end{lemma}

\begin{proof}
The first part follows immediately from the definition of recoding.
To prove the second statement,
note that if $\alpha \in \cA'$,  then $\alpha$
determines a coloring of $\{(0,0)\}$ for the shift $Y$.  Since
$\{(0,0)\}$ $\Psi^{-1}$-codes $F$,  it follows that $\alpha$ determines a unique coloring
$\beta$ of $F$ for $X$.
 The assignment $\alpha \mapsto \beta$ determines
a bijection  from $\cA'$ to  $\cA_{F}$ which, as a block map, determines
an isomorphism $\Psi_F \circ \Psi^{-1}$ from $(Y, \cA')$ to $(X_F, \cA_F)$.
\end{proof}

It is frequently useful to know that a finite set coded by a set $A$
is also coded by a finite subset of $A$.  This follows via an easy compactness
argument:
\begin{lemma}\label{lem: compactness}
Assume that $X$ is a $\Z^d$-subshift.  If $A\subset \Z^d$ codes $B$ and $B$ is finite,  then there is a finite subset $A_0 \subset A$
such that $A_0$ codes $B$.
\end{lemma}

\begin{proof}
Without loss of generality, it suffices to prove the result when $B$
contains a single point $ b \in \Z^d$ which is not  an element of $A$.  If
the result fails, then for every $m \ge 0$ there exist
$x_m, y_m \in X$ such that $x_m( b) \ne y_m( b)$, but $x_m( u) = y_m( u)$
for all $ u \in A$ with $\| u \| \le m$.  Since $X$ is compact, by passing to by
subsequences if needed, we can assume that $\lim_{m\to\infty} x_m = x'$ and $\lim_{m\to\infty} y_m = y'$ for some
$x',y'\in X$.  Then $x'( b) \ne y'( b)$, but
$x'( u) = y'( u)$ for all $ u \in A$, a contradiction as
$A$ codes $\{ b\}$.
\end{proof}

\subsection{Notions of parallel}
\label{sec:parallelisms}
We summarize the various notions of parallel that we use throughout the sequel.

By a  {\em ray} in $\R^2$, we mean 
a translate of the set $\ell_v = \{tv \colon t \ge 0, \ v\ne 0 \in \R^2\}$, 
and we view a ray as oriented in the direction of increasing $t$.

Two rays $\ell_1$ and $\ell_2$ are {\em positively parallel} if one is a translate of the other, and
they are {\em antiparallel} if $\ell_1$ and $-\ell_2$ are positively parallel.  
We say that two rays are
{\em parallel} if they are either positively parallel or antiparallel.  

We extend these conventions to line segments, and 
we say that an oriented line segment $J$ is
{\em positively parallel } to a ray $\ell$ if a translate of $J$ lies in $\ell$ with matching
orientations, and we say that the orientated line segment $J$ is {\em antiparallel}  if $J$ is positively parallel to $-\ell$.  

Similarly, we say that two oriented
line segments are {\em positively parallel} if a translate of one lies in the other
with matching orientations and are {\em antiparallel} if one is positively parallel with  the
other with reversed orientation.  Since we need to distinguish the various notions our terminology differs a bit  from that in~\cite{CK}, where parallel 
corresponds to our use of positively parallel, while the use of antiparallel is the same.

A polygon $\cP\subset \R^2$ inherits an orientation from $\R^2$, and this orientation induces an
orientation on the boundary $\partial \cP$, and this further restricts to an orientation on
each edge of $\cP$.  
For a convex polygon, no two edges can be positively parallel, but pairs of edges may
be antiparallel. 

\subsection{Expansive and nonexpansive}
The fundamental concept related to one region coding another is
that of expansivity, defined in Milnor~\cite{M} and developed by~\cite{BL}, and we review this in our
particular setting of two dimensions.  Letting $d$ denote the distance in $\R^2$,
a subspace $L$ of $\R^2$
is {\em expansive} if there exists $r>0$ such that the
$r$-neighborhood $N_r = \{u\in \R^2\colon d( u, L) < r\}$ of $L$
codes $\R^2$
(the analogous definition can be made in any dimension).
Any subspace that is not expansive is called a {\em nonexpansive} subspace.

Nonexpansive subspaces are common:
\begin{theorem}[Boyle and Lind~\cite{BL}]
\label{th:BL}
If $X$ is an infinite compact metric space with a continuous $\Z^k$-action, then for each
$0\leq j < k$ there exists a $j$-dimensional subspace of $\R^k$ that is nonexpansive.
\end{theorem}
For the two dimensional setting, an immediately corollary is that a system $X$ is
finite (and hence doubly periodic)
if and only if every subspace of $\R^2$ is expansive.

For our purposes, the notion of expansiveness can be refined to consider one-sided
expansiveness, where the coloring of $N_r$ determines the coloring of
one component of the  complement of $L$.  We make this more
precise (similar notions were considered in~\cite{BlM, BM, CK}):
\begin{lemma}
\label{lemma-half}
Assume $X$ is a $\Z^2$-subshift
and suppose $H$ is an (open or closed) half space in $\R^2$.
Then either $H$ codes all of $\R^2$ or $H$ codes itself
but no points of $\Z^2 \setminus H$.
In particular, if any subset of $H$ codes any point of $\Z^2 \setminus H$,
then $H$ codes all of $\R^2$.
\end{lemma}

\begin{proof}
Suppose there is no  $b \in \Z^2 \setminus H$ such
that $H$ codes $\{b\}$.  In this case, $H$ codes subsets of itself
and no other subsets of $\R^2 \cap \Z^2$.
Otherwise,  there exists $b \in \Z^2 \setminus H$
such that $H$ codes $\{b\}$.
We prove this implies
$H$ codes $\R^2$.

First consider a special case: assume  that $H$ is closed
and there exists some $z \in \partial H \cap \Z^2$.
Let $w = b -z \in \Z^2$ and define $H_1 = w + H$.  Then  $\partial H_1
= L+w$ contains $b$.  Thus $H_1$ is a closed half space properly
containing $H$. We claim that if $u \in H_1 \cap \Z^2$,
then $H$ codes $\{u\}$.  To see this, let $v = u -b$.  Since $b \in \partial H_1$
and $u \in H_1 \cap \Z^2$, we have $v + H_1 \subset H_1$ and hence
$v + H \subset H$. Since $H$ codes $\{b\}$, we have that
$v +H$ codes $v +b = u$. But $v + H \subset H$, proving the claim.

By the claim, it follows that  $H$ codes $H_1$.
Define $H_n =  nw+H$.  Then since $H$ codes $H_1$, it follows from the translation
invariance (Remark~\ref{rem:translate}) that  $H_n = H+nw$ codes $H_{n+1}  = H_1 +nw$.
Hence $H$ codes $\bigcup_n H_n = \R^2,$   meaning that $H$ codes $\R^2$.
This completes
the proof in the special case that $H$ is closed and $\partial H \cap Z^2$ is nonempty.

We now turn to the general case, assuming that $H$ is an open or closed
half space bounded by $L$ (and no assumption that the line $L$ contains
points of $\Z^2$).
The point $b \in \Z^2$ is coded by $H$, but $b \notin H$.
By Lemma~\ref{lem: compactness}, there is a finite set $A \subset H$ which
codes $\{b\}$.
For each $a \in A$,  let $Y_a$ denote the closed
half space contained in $H$ whose boundary is the line $L_a$ which is parallel to
$L$ and contains $a$.  If
\[
Y = \bigcup_{a \in A} Y_a,
\]
then $Y \subset H$ and $Y$ is a closed half space which codes $\{b\}$ and
$b \notin Y$.  Also $\partial Y$ contains some point of $A$ and hence
some point of $\Z^2$.   It follows that $Y$ satisfies the hypothesis
of the first  case and so $Y$ codes $\R^2$.  Since $Y \subset H$,
we also have that $H$ codes $\R^2$.
\end{proof}

Note that for any $v \in \R^2$ (not necessarily integral),
a half space $H$ is expansive if and only if $v + H$ is expansive.
For $v \in \Z^2$, this follows immediately from the translation invariance (Remark~\ref{rem:translate}).
More generally, for any $v \in \R^2$ there exist $z_1, z_2 \in \Z^2$
such that $z_1 + H \subset v+ H \subset z_2 + H$ and so $z_1 + H$
expansive implies $v + H$ is expansive and
$v + H$ expansive implies $z_2 + H$ is expansive.

We encapsulate the  dichotomy of Lemma~\ref{lemma-half} in the
following definition:

\begin{defn}\label{def expansive}
Assume $X$ is a $\Z^2$-subshift.
If $H$ is an (open or closed) half space in $\R^2$,  we  say that
$H$ is {\em expansive} if $H$ codes $\R^2$ and otherwise we say that
$H$ is {\em nonexpansive.}
If $H$ is expansive and  $\ell$ is a ray parallel to the boundary of $H$ whose
orientation agrees with the orientation $\partial H$ inherits from
the standard orientation on $H$, 
we say that $\ell$ is an {\em expansive ray in $X$} and otherwise we say that
$\ell$ is a {\em nonexpansive ray in $X$.} When it is clear from the context, 
we shorten this and say that $\ell$ is {\em expansive} (or {\em nonexpansive}). 
\end{defn}

\begin{remark}\label{rem:iso invariance}
It is easy to see that if $(X, \cA)$ and $(Y, \cA')$ are isomorphic
shifts, then a ray $\ell$ is expansive for one if and only if it is
expansive for the other (see Remark~\ref{rem:translate}).
\end{remark}

We note that the half space $H$ being nonexpansive
is equivalent to the existence of $x_1, x_2 \in X$  with $x_1 \ne x_2$ such that
$x_1(i,j) = x_2(i,j) $ for all $(i,j) \in H$.  This non-uniqueness in the extension of the half space
is often how we make use of this notion.

A one-dimensional subspace $L$ of $\R^2$ is (two-sided) expansive if
for some $r >0$, the strip $N_r(L) = \{ u\in \R^2\colon d( u, L) \le r\}$
codes $\R^2.$ This implies that the action on $X$ by any nonzero
element $v \in L \cap \Z^2$ is an expansive homeomorphism of $X$.
The following corollary shows that a subspace $L$ is expansive in this
sense if and only if both of its complementary half spaces satisfy
our definition of one-sided expansiveness
(Definition~\ref{def  expansive}).

\begin{prop}
If $H$ and $H'$ are the two closed half spaces whose
common boundary is $L$ (so $H \cup H' = \R^2$) and $H$ codes all of $\R^2$,
then there exists $r >0$ such that the closed strip
$N_r(H)  = \{u \in H \colon d(u, L) \le r\}$ codes
all of $H'$.
\end{prop}

\begin{proof}
Without loss of generality, we can assume that
$L = \partial H$ contains some point of $\Z^2$: choosing a
(not necessarily integral) translate $L_0$ of $L$ which
lies in $H$ and does contain a point of $\Z^2$, we can prove the result
for $H_0 \subset H$ with $L_0 = \partial H_0$ and
obtain the result for $H$ (possibly with a larger value of $r$).

As in the special case in the proof of Lemma~\ref{lemma-half},
$H$ codes  some $b \in \Int(H') \cap \Z^2$.
Choose $z \in L \cap \Z^2$ and set $w = b-z$. There is a finite set
$A \subset H$ which codes $\{b\}$.  Let $\delta = d(b, L)$. Suppose
$u \in H'$ and $d(u, L) \le \delta$.  Setting $v = u-b$, the component
of $v$ orthogonal to $L$ has length $\le \delta$ and so
$A + v \subset N_r(H)$ where $r = \delta + \diam(A)$.  Also $A + v$
codes $b +v =u$ and so $N_r(H)$ codes the closed strip $S$ whose boundary
components are $L$ and $L +w$.  The strip $S$ is parallel to $L$
and has width $\delta$.  The same argument shows that
$S \cup N_r(H)$ codes the strip $S +w$.  Inductively,  it follows that
$N_r(H) \cup (nw +S)$ codes $N_r(H) \cup ((n+1)w +S)$, and so
$N_r(H)$ codes $H'$.
\end{proof}

The set of expansive rays in $\R^2$ is open (see~\cite{BL, CFK}). This also
follows immediately from Lemma~\ref{lem: compactness}, which
gives the existence of a finite set $A \subset H$ which codes
$b \notin H$, and the fact that the set of oriented rays in the plane which
span lines separating $b$ from $A$ is an open set.

It thus follows that the set of  nonexpansive rays is closed,  and it is known to be nonempty
if $X$ is infinite (see Theorem~\ref{th:BL}).
For the full shift $\cA^{\Z^2}$, it is easy to see that all rays are nonexpansive;
there are no expansive half spaces.  The nonexpansive rays
play a significant role in the dynamics of $\Z^2$-subshifts
because the boundary of a nonexpansive half space  forms a barrier to coding.
In particular, Lemma~\ref{lemma-half} asserts that if $H$ is nonexpansive,
then no subset of $H$ can code a subset of
$\Z^2 \setminus H$.

\section{Defining the class of shifts}
\subsection{Polygonal shifts}
\label{sec:polygonal}

We have assembled the tools to define the class we study: 
\begin{defn}\label{def polygonal}
Suppose  $X$ is an infinite  $\Z^2$-subshift,  $\cP$ is a convex  integer polygon, 
and $v$ is a vertex of $\cP$.  If
$\cP \setminus \{v\}$ $X$-codes $\{v\}$, 
then we say that $\cP$ is a {\em coding polygon} for the vertex $v$.
If $\cP$ is coding for each of its vertices, 
we say $X$ is  {\em polygonal} with respect to $\cP$ or
that $\cP$ is a {\em coding polygon} for $X$.

A polygonal $\Z^2$-system is  {\em triangular} if the associated polygon is
a triangle.
\end{defn}

Note that translation invariance implies 
that when $\cP \setminus \{v\}$ $X$-codes $\{v\}$, we also have that 
$(\cP + u) \setminus \{v+ u\}$ $X$-codes $\{v +u\}$
for all $u \in \Z^2$.  Thus it makes sense to discuss a coding polygon defined
only up to translation in $\Z^2$.
However, coding polygons, even up to this translation, are not unique and
it takes work to understand to what extent a coding polygon can be simplified. 
One notion of simplification is having the fewest number of edges, 
and this motivates us to restrict our attention to
convex polygons.  If a non-convex polygon is coding, then
its convex hull has fewer sides and is also a coding polygon.

\begin{prop}\label{prop: polygon sides}
Suppose  $X$ is a  $\Z^2$-subshift and  $P$ is a coding  polygon for $X$.
If $\ell$ is a nonexpansive
ray in $X$, then $\ell$ is positively parallel to an edge of $P$ whose orientation matches
the orientation of $\ell$.
\end{prop}

\begin{proof}
Let $L$ be the one-dimensional subspace of $\R^2$ containing $\ell$
and let $H$ be the open half space bounded by $L$ such that 
expansiveness of $H$ implies expansiveness of $\ell$.
Suppose first that $L$ is not parallel to any edge of $P$.
Then there is vertex $e$ of $P$
such that $P \cap (e +L) = \{e\}$ and $P \subset e +H$. 
Since $P \setminus \{e\} \subset e+H$ codes $e \notin e+H$, 
Lemma~\ref{lemma-half} implies that $e +H$ is expansive.
It follows that if $\ell$ is nonexpansive, then it is parallel
to an edge of $P$; we are left with showing that there is an edge which is positively parallel to $\ell$.
If it is positively parallel to one edge and antiparallel to another, then those
edges have opposite orientations and so $\ell$ is positively parallel to
one of them.  Finally, if $\ell$ is antiparallel to a single edge $E$, 
then there is a unique vertex $e \in P$ and a
unique supporting line $L$ parallel to $\ell$ such that $L \cap P = \{e\}$.
If $H$ is the open half space which is bounded by $L$ 
and which  contains $P \setminus \{e\}$,  
Lemma~\ref{lemma-half} implies that $H$ is expansive
(note that  $P \setminus \{e\} \subset H$ codes $e \notin H$).
The orientation $L$ inherits
from $H$ is the opposite of the orientation $E$ inherits from $P$.
Since  the ray $\ell$ is antiparallel to $E$ and $H$ is expansive, 
$\ell$ must be expansive, a contradiction.  The only remaining
possibility is that $\ell$ is positively parallel to $E$.
\end{proof}

Although coding polygons are not unique, the existence of a coding
polygon implies that scaled versions are also coding polygons.  To make this precise, given $P\subset\Z^2$, we
write
$$nP = \{nx\colon x\in P\}.
$$ 
We frequently make use of the following straightforward
observation: 

\begin{obs}\label{scale up}
If $X$ is polygonal with
respect to the convex polygon $\cP$, then 
it is also polygonal with respect to the polygon $n \cP$ for every $n \in \N$. 
\end{obs}

This can be seen  by noting that if $v$ is a vertex of $\cP$, 
there is a translation $T$ such that $T(v) = nv$ and
then $T(\cP)$ is a subpolygon of $n\cP$ whose vertex at $nv$ coincides with that
of $n\cP$.

\subsection{Examples of polygonal shifts}
\label{sec:examples}
We give various examples of polygonal shifts.  
\begin{example}
\label{ex:ledrap}
{\em Ledrappier three-dot system}~\cite{ledrappier}.
Let $\CA = \{0,1\}$ be the field with two elements
and take $X$ to be the subshift of $\CA^{\Z^2}$ 
defined by 
\[
X = \{x \in \CA^{\Z^2} \colon  x(i,j) + x(i+1, j) + x(i, j+1) = 0 \mod 2\}
\]
for $i, j\in\Z$.  
Note that if $x \in X$ and
$R_i(x)$ is the element in the one-dimension shift
$\cA^Z$  obtained by restricting $x$
to its $i^{th}$ horizontal row, then $R_{i+1}(x) = \phi(R_i(x))$
where  $\phi$ is the endomorphism defined by
$\phi(y)_0 = y_0 + y_1 \pmod 2$.  
The $\Z^2$-subshift $X$ is triangular (with respect to the triangle $T$ 
with vertices $(0,0), (1,0),$ and $(0,1)$).
It has three nonexpansive rays, 
which are the positive $x$-axis, the negative $y$-axis, 
and the ray $(-t, t),\ t \ge 0$.
\end{example}

Ledrappier's three dot system and related
algebraic systems have been studied by Ledrappier~\cite{ledrappier},
Einsiedler~\cite{E}, and Kitchens and Schmidt~\cite{KS}.  In particular we have the following extension: 
\begin{example}
\label{ex:einsiedler}
{\em Einsiedler's examples.} 
In~\cite{E} Einsiedler proves the existence of closed $\Z^2$-invariant subsystems
of the Ledrappier example
realizing any horizontal directional entropy between $0$ and $\ln(2)$. 
Since these are subsystems of the Ledrappier system, they are all triangular with respect
to the triangle $T$.  Following~\cite{E} and~\cite{KS2}, 
we describe one such example.
Taking $X$ to be the Ledrappier system of Example~\ref{ex:ledrap}, 
consider the subset 
$$Y_0 = \{ x \in X \colon x(2v) = 0 \text{ for all } v \in \Z^2\}.
$$
Then $Y_0$ is invariant under the 
$\Z^2$ -action
obtained by restricting the
$\Z^2$-action on $X$ to the lattice $(2\Z)^2$.  
While $Y_0$ is not invariant under the full $\Z^2$-action,
defining $Y_1 = Y_0 + (1,0),\  Y_2 = Y_0 + (0,1), \ Y_3  = Y_0 + (1,1),$
then  $Y = Y_0 \cup Y_1\cup Y_2\cup Y_3$ 
is a closed proper $\Z^2$-invariant subset of $X$.
Let $R_i$ denote the restriction of $Y_i$ to the  $x$-axis and $R$ denote the restriction
of $Y$. Then each $R_i$ is a closed subset of the
 one-dimensional full shift space $\Sigma = \cA^{\Z}$ with
\begin{align*}
R_0 &= \{y \in \Sigma \colon y_{2n} = 0 \text{ for all } n \in \Z\}\\
R_1 &= \{y \in \Sigma \colon y_{2n+1} = 0 \text{ for all } n \in \Z\}\\
R_2 &= \{y \in \Sigma \colon y_{2n} = y_{2n+1} \text{ for all } n \in \Z\}\\
R_3 &= \{y \in \Sigma \colon y_{2n} = y_{2n-1} \text{ for all } n \in \Z\}.  
\end{align*}
Define $\sigma\colon R \to R$ to be the left shift and  observe that $\sigma^2(R_i) = R_i$.
One checks easily that $\sigma^2|_{R_i}\colon R_i \to R_i$ is 
conjugate to the full $2$-shift, and so $\sigma^2\colon R_i \to R_i$ has topological entropy
$\ln(2)$.  Since $R_i \cap R_j$ contains at most the two
$\sigma$-fixed points $\bar 0$ and $\bar 1$ for all $i\ne j$, it follows that 
$\sigma^2\colon R \to R$ has topological entropy $\ln(2)$ and hence
$h(\sigma) = \ln(2) / 2$.
\end{example}

It follows from Einsiedler's results that uncountably many horizontal
entropies can be realized in constructing the examples in~\ref{ex:einsiedler}.  
All but countably many of the associated
subshifts are not sofic, since there are at most countably many
sofic systems with a given alphabet.  Thus some of the subshifts realized in Example~\ref{ex:einsiedler} are not
sofic and, in particular, are not subshifts of finite type.

\begin{example} {\em Low complexity examples.}
Recall that, by convention, the alphabet $\cA$ only contains letters which occur in the
language of $X$.
Polygonal systems arise naturally in studying the Nivat Conjecture, 
and in this direction, 
it follows immediately from~\cite[Corollary 2.6]{CK} that (note our terminology differs, and related results appear in~\cite{CG, KM}): 
\begin{prop}\label{prop: nivat}
Suppose $X$ is a $\Z^2$-subshift with alphabet $\cA$,
$S$ is a finite convex subset of $\Z^2$, and $C(S)$ denotes the
number of legal $X$ colorings of $S$. If
$$C(S) \le |S| + |\cA| -2,$$
then $X$ is polygonal with a coding polygon  which can be chosen
as a subset of $S$.
\end{prop}
In particular, it follows from~\cite{CK} that any counterexample to the Nivat conjecture must be polygonal.
\end{example}

\begin{example} {\em Non-abelian groups.}
Similar to the construction of the Ledrappier system,  
one can take a finite (possibly non-abelian) group $G$
as the alphabet and require, for example, that the product of the
colors at the vertices of a convex polygon $P$ is the identity (or
some other fixed $g \in G$). 
\end{example}

\begin{ex}\label{enveloping SFT}{\em Shifts of finite type.}
We claim that any polygonal shift $X$ can be written 
as a countable intersection of shifts of finite type, each of which is polygonal with the same polygon as $X$.  

To check this, first note that any $\Z^2$-subshift $X$ can be written as a countable intersection of shifts of finite type: 
namely, $X$ can be defined as all colorings which do not contain any
elements of a countable set $E$ of excluded block colorings (excluding only 
 finitely many elements of $E$ results in a shift of finite type).  
 By excluding larger and larger finite subsets $E_n$ of $E$, we obtain a nested sequence $X_n$ of
shifts of finite type, each of which contains $X$.  If the sets $E_n$ are chosen such that
$E = \bigcup_{n\in\N} E_n$, then the  intersection of
the resulting shift of finite type $\bigcap_{n\in\N} X_n$ is $X$.

Now suppose $X$ is polygonal with coding polygon $\cP$ and $F$ is the finite set of
colorings of $\cP$ which do not occur in $X$ (and so excluding $F$
  incorporates the subshift of finite type constraints given by the fact that $\cP$ is
  a coding polygon). Suppose $E$ is the countable
set of excluded block colorings defining $X$. 
By choosing $E_n \subset E$, $n\ge 1$
such that $F \subset E_n$ and $E = \bigcup_{n\in\N} E_n$, 
then each of the shifts of finite type  $X_n$, defined by excluding the blocks
$E_n$, is polygonal with coding polygon $\cP$.  Thus the polygonal shift $X$ with coding polygon $\cP$
is a countable intersection of shifts of finite type $X_n$, 
each of which is polygonal with coding polygon $\cP$.
\end{ex}

\begin{example}{\em Products.}
If $X_1$ and $X_2$ are  $\Z^2$-shifts with alphabets $\cA_1$ and
$\cA_2$, then their Cartesian product is the $\Z^2$-shift $Y$ with alphabet
$\cA_Y := \cA_1 \times \cA_2$ consisting of $y$ such that $p_1(y) \in X_1$ and
$p_2(y) \in X_2$ where each $p_i\colon Y \to X_i$ is the map induced by projecting
$\cA_1 \times \cA_2$ onto the $i^{th}$ component.
\end{example}

If $X_1$ and $X_2$ are polygonal $\Z^2$-shifts with the same coding polygon
$P$, then it is immediate that $X_1 \times X_2$ is also polygonal with
coding polygon $P$.  In other words the polygonal shifts with a fixed coding
polygon form a semi-group under Cartesian product.
More generally we have: 
\begin{proposition}\label{prop: cartesian}
If $X_1$ and $X_2$ are polygonal $\Z^2$-shifts
with respect to polygons $P_1$ and $P_2$,
then $X_1 \times X_2$ is also polygonal. 
\end{proposition}
\begin{proof}
Consider the positively oriented edges $\{v_i\}$ of the two polygons $P_1$ and $P_2$ as vectors, 
 and order them that they form the edges of a convex polygon $P$; more precisely, order
these edges in the circular order determined by the angles they form with the
$x$-axis determines a convex polygon.  Then let
$e_i$ denote the segment from $\sum_{j \le i}v_j$ to
$\sum_{j \le i+1}v_j$, concatenating these edges in order gives a curve
and this curve is closed because the sum of the edges in
each of $P_1$ and $P_2$ is zero.  Because of the order the curve is
the boundary of a convex polygon. Each edge $e_i$ is positively parallel to
to the vector $v_i$.

If there are positively parallel edges, one  in $P_1$ and the other in  $P_2$, 
this creates successive parallel edges in the new polygon $P$; to simplify, 
we delete the vertex between
them to create single edge of $P$ whose length is the sum of
the lengths of the two parallel edges.

Recall the alphabet for $X_1 \times X_2$ consists of ordered
pairs of colors from the alphabets of $X_1$ and $X_2.$
To check that $X_1 \times X_2$ is a polygonal system, consider a
coloring of all but one vertex $w$ of $P$.  
We claim that we can translate the
polygon $P_1$ associated to $X_1$  such that $P_1$ lies in $P$ and
a vertex of $P_1$ coincides with $w$ (and the analogous statement holds for $P_2$).
The coloring of this copy of $P_1$  (for system $X_1$)
with $w$ deleted uniquely determines the color
of the vertex $w$ of $P_1$ and hence the first component of the pair which
is the coloring for $X_1 \times X_2$.  The second component of the
coloring for $w \in P$  is obtained similarly.

We are left with proving the claim, showing that a translate of $P_1$ lies in $P$ with
a vertex at $w$. 
Without loss of generality, we can assume that no edge of $P$ is horizontal,  the
vertex $w$ of  $P$ is at the origin, and the remainder of $P$
lies in the upper half plane. We also assume that our ordering of
the $\{v_i\}$ starts with the edge emanating from $0$ in the positive orientation of
the edges of $P$. Then there is a translate of $P_1$
with a vertex at $0$ and such that the edges of $P_1$ incident to $0$ lie in $P$ (otherwise
we have contradicted the ordering on the edges of $P$).   
Let $q$ be the highest vertex (meaning in the $y$ direction) 
of $P$.  Then  $q = \sum_{i = 0}^k v_i$ where
$v_i$ has a positive $y$ coordinate for $0 \le i \le k$ and a negative
$y$ coordinate for $i> k$. Thus $0$ and $q$ divide the boundary of $P$ into
two pieces: $K^+$ where all the oriented edges have a positive $y$ component and $K^-$
where they all have a negative $y$ component.

Denote the edges of $P_1$ by $\{u_n\}$. 
By construction, the beginning of $u_0$,
the  first edge of $P_1$ starting at $0$, lies in $P$ or on its boundary.
Suppose now that the boundary of $P_1$ intersects
$K^+$ at a point $z$, and then crosses out of $P$.
In other words, suppose $u_{i_0}$ and $v_{j_0}$ are edges of $P_1$ and $P$ respectively,
  containing the point $z$ and with $u_{i_0} \ne v_{j_0}$.
  If $z$ is a vertex of $P_1$  (respectively $P$), we choose
$u_{i_0}$  (respectively  $v_{j_0}$) such that $z$ is the beginning endpoint of
$u_{i_0}$  (respectively $v_{j_0}$).
Since $u_{i_0}$ is exiting the polygon $P$, the angle with respect to the $x$-axis is greater for  $v_{j_0}$
than for $u_{i_0}$, 
Then
\[
z =  a v_{j_0} + \sum_{i = 0}^{j_0 -1} v_i
\]  for some $0 \le a , 1$.
But also since
the edges of $P_1$ are $\{u_n\}$
\[
z =  b u_{i_0} + \sum_{i = 0}^{i_0 -1} u_i
\]  for some $0 \le b < 1$.
Note  that because of the ordering of edges, every $u_n$ for $0 \le n \le i_0$
must be equal to some $v_k$ with $0 \le k \le j_0$.
In particular $u_{i_0} = v_{k_0}$ for some $0 \le k_0 < j_0$.
It follows that either
every edge $u_n$, $n \le i_0$, is an edge of both $P$ and $P_1$ or the $y$-component
of $a v_{j_0} + \sum_{i = 0}^{j_0 -1} v_i $ 
is strictly greater than  the $y$-component of  $b u_{i_0} + \sum_{i = 0}^{i_0 -1} u_i $
The second condition can not hold, since both sums equal $z$.
So if the boundary of $P_1$ intersects $K^+$ in more than the vertex
$0$, it can only do so in an arc of edges
which are common to $P_1$ and $P$.  A similar argument applied to $K^-$ shows
that if the boundary of $P_1$ intersects the boundary of $P$ it
can only do so in an arc of edges
which are common to $P_1$ and $P$.  This implies $P_1 \subset P$.
The same argument shows $P_2 \subset P$.
\end{proof}

\begin{prop}\label{scale down}
If $X$ is polygonal with
respect to the convex integer polygon $\cP$ and $\cP = n \cP_0$ for some integer polygon 
$\cP_0$ and integer $n >1$,
then there is a recoding $Y$ of $X$ which is polygonal with respect to
$\cP_0$.  Hence $m \cP_0 = \frac{m}{n}\cP$ is also a coding polygon for $Y$.
\end{prop}

This result shows that if $\cP$ is a coding polygon for $X$ and
$\cP = n \cP_0$ for som $n >1$, then  $X$ can be recoded to $Y$ with a strictly
smaller but similar coding polygon which is primitive (meaning that it is 
not an integer multiple of a smaller integer polygon).

\begin{proof}
Without loss of generality we may assume $0$ is a vertex of $\cP_0$ (and hence also of $\cP$). 
Let $(e_i)_{i=0}^k$ denote the edge vectors of $\cP$ starting at $0$
and taken in a counter-clockwise order.

Then  $ v_j = \sum_{i=0}^j e_i$ is the $j^{th}$ vertex of $\cP$ in this ordering, and   
$v_k = \sum_{i=0}^k e_i = 0.$
Since $\frac{1}{n}\cP$ is an integer polygon, so is
$\frac{m}{n}\cP = m \cP_0$ for $1 \le m \le n$.  
In particular, each edge of
$ m \cP_0$ is a segment with endpoints in $\Z^2$.

For  a line segment $e$ in $\R^2$ whose endpoints lie in $\Z^2$, set $\mu(e) = |e\cap\Z^2|-1$, 
where $|\cdot|$ denotes the number of points. 
This particular choice of definition for $\mu$ is taken such that we have 
$\mu(n e) = |n|\mu(e)$ for all $n \in \Z$ and if
$e$ and $f$ are line segments intersecting only in a common endpoint, 
then $\mu( e \cup f) = \mu(e) + \mu(f)$.  Note that if $w \in \Z^2$, 
then $\mu(e +w) = \mu(e)$.

Set $\eps_i  := \frac{1}{n} e_i $, meaning that $\eps_i$ is  
the $i^{th}$ edge of $\cP_0$.  Then $\mu(\eps_i) = 
\frac{1}{n} \mu(e_i)$ .
Define  $\cP_1 :=\frac{n-1}{n}\cP = (n-1) \cP_0$.
Let  $\eta_i := (n-1) \eps_i$  denote the  $i^{th}$ edge of $\cP_1$
and so $\mu(\eta_i) = (n-1) \mu(\eps_i)$.
Then for each $i$, we have 
$\mu(\eta_i) + \mu(\eps_i) = (n-1) \mu(\eps_i) + \mu(\eps_i) = \mu(e_i)$.
Thus, for each edge $e_i$ of $\cP$, 
there are exactly $\mu(\eps_i)$ translates in $\Z^2$ 
of $\cP_1$ each of which lies in $\cP$ and has an integer translate of
the edge $\eta_i$ of $\cP_1$  lying in $e_i$. 

Let $\Psi\colon X \to X_{\cP_1}$ be the canonical recoding of $X$ (see Definition~\ref{def: recoding}) 
via $\cP_1$ and let $Y = X_{\cP_1}$. Then the polygon $\cP$ \ $\Psi$-codes
a translate of $\cP_0 = \frac{1}{n} \cP$.  Likewise $\cP_0$ \ $\Psi^{-1}$-codes a
translate of $\cP$.  It follows that $\cP_0$ is a coding polygon for
$Y$.  Hence by Observation~\ref{scale up}, 
$m \cP_0 = \frac{m}{n}\cP$ is also a coding polygon for $Y$.
\end{proof}

\subsection{Refining notions of expansivity}

Suppose $L$ is a rational line in $\R^2$ containing a point of $\Z^2$ (and hence infinitely
many points of $\Z^2$).  Recall that $L + \Z^2$ is a discrete set of lines, meaning there exists 
$r >0$ such that any line $z + L$ distinct from $L$ and with $z \in \Z^2$ must have
distance from $L$ equal to $mr$ with $m \in \N$. There are two closest
integer translates of $L$ which have distance 
$r$ from $L$, lying on opposite sides of $L$.

If $L\subset\R^2$ is a one dimensional subspace, we refer to the intersection
of a connected segment of $L$ with $\Z^2$ as a {\em block} in
$L \cap \Z^2$.

\begin{defn}\label{def closing}
Assume $L\subset\R^2$ is a one dimensional subspace with rational slope and suppose
$L$  bounds a nonexpansive closed half space $H$.  Let $L_0$ 
be the closest line of the form $z + L$ in the complement of $H$ for some $z \in \Z^2$.
If there exists $N>0$ such that every  block  $\B$ in $L_0 \cap \Z^2$ of length $\ge N$
the set $H \cup \B$ codes $H \cup L_0$  then, we say that 
$H$ is {\em  closing.}  If $\ell$ is the ray in $L$ whose orientation is inherited
from $H$, we say that $\ell$ is {\em closing}.
\end{defn}

Note that by definition, a ray that is closing is also
nonexpansive and has rational slope.
To explain the rationale behind the use of the term closing, note, for
example, that in the Ledrappier system (Example~\ref{ex:ledrap}), the
upper half space $H = \{(u,v) \colon v \ge 0\}$ is nonexpansive.  The
subspace $H$ is also closing.  This latter property is equivalent to
the fact that the endomorphism $\phi$ defining the system is both
right and left closing in the sense of~\cite[Chapter 8]{LM}.

We note the relevance of the property of closing to  polygonal shifts:
\begin{prop}\label{prop: closing}
If $X$ is polygonal with coding polygon $\cP$ and $\ell$ is a
nonexpansive ray with the same direction as an oriented edge
of $\cP$, then $\ell$ is closing.
\end{prop}
\begin{proof}
Without loss of generality, we can assume that the ray $\ell$ is the positive
horizontal axis, meaning that the oriented edge lies in the horizontal axis $L$ and
$\cP$ lies in the closed upper half space $H$
with the oriented edge of $\cP$ lying in $L$ matching the orientation of $L$.
Then $L_0$ is the
line $L + (0,-1)$.  Let $N$ be the number of points in $J := L \cap \cP$.
Set $\fB = J (0,-1)$ and note that every integer point of
the polygon $\cP + (1,-1)$ lies in $H \cup \fB$ except one, namely the first
point $b$ to the right of  $\fB$ in  $L_0$.  Since $\cP$ is
a coding polygon, the coloring at $b$ is determined by
the coloring of $(\cP + (1,-1)) \setminus \{b\}$
and hence by $H \cup \fB$. Repeating this, 
it follows that $H \cup \fB$ codes all points
to the right of $\fB$.  A similar argument shows it codes all points to
the left of $\fB$.
\end{proof}

It follows from Remark~\ref{rem:iso invariance} that if a ray is nonexpansive for $(X, \cA)$, then it is
also nonexpansive for any isomorphic $\Z^2$-shift $(Y, \cA')$.  Our next lemma
shows that a recoding (and its inverse) preserves closing rays: 
\begin{lemma} \label{lem:recoding expansive}
Suppose $\Psi\colon (X, \cA) \to (Y, \cA')$ is a recoding via a finite set $F$
and suppose $\ell$ is a rational nonexpansive ray in $\R^2$.
Then $\ell$ is closing for $X$ if and only if it is 
closing for $Y.$
\end{lemma}
\begin{proof}
Suppose $\ell$ is closing for one of $X$ or $Y$. We show it is closing for
the other.  By a change of coordinates, without loss of generality we may 
assume that $\ell$ is the
positive $x$-axis. Let $L$ be the $x$-axis and let $H$ be the closed
upper half space with boundary $L$.  Thus by our
  hypothesis, $H$ is nonexpansive for $X$.

By Lemma~\ref{lem:recoding}, if $T$ is an action on $X$ induced by translating
by some element of $\Z^2$, then recoding via $T(F)$ is the same as
recoding via $F$ and then translating by $T$.
Since translating by $T$ preserves closing
half spaces, we can assume that $F$ lies in $H$ and contains $(0,0) \in L$, but
contains no point of $\Z^2 \setminus H$.  

Let $L_0 = L +(0,-1)$ and
suppose  $B$ is  a finite block in $L_0$.  Since $F \subset H$ and
$F$ \ $\Psi$-codes $\{(0,0)\}$, by translating  it follows that $H$ \ $\Psi$-codes $H$.  
Thus an $X$-coloring of $H$ determines a $Y$-coloring of $H$.
Likewise  $H +(0,-1)$ \ $\Psi$-codes $H+(0,-1)$. By the definition of  recoding, 
$\{(0,0)\}\ \Psi^{-1}$-codes $F$ and therefore
$\{(i,j)\} \ \Psi^{-1}$-codes $F+(i,j)$.  
Thus it follows that for all $i,j$, we have 
that $(i,j)+F\ \Psi$-codes $\{(i,j)\}$ and
$\{(i,j)\}\ \Psi^{-1}$-codes $F+(i,j)$ and hence code $\{(i,j)\}$.  

Suppose now that $H$ is $X$-closing and $B \subset L_0$ is a block
such that $H \cup B$ \ $X$-codes $L_0$.  We claim that 
$H \cup B$ \ $Y$-codes $L_0$,  
Since $\{(i,j)\}\ \Psi^{-1}$-codes  $\{(i,j)\}$, we have that $H \cup B$ \ 
$\Psi^{-1}$-codes $H \cup B$ which $X$-codes $H \cup L_0 = H +(0,-1)$.
This in turn $\Psi$-codes $H +(0,-1) \supset L_0$ and we have
shown that $H \cup B$ \ $Y$-codes $L_0$. 
This proves the claim and it follows that $H$ is $Y$-closing.

Conversely, suppose $H$ is $Y$-closing
and let  $B \subset L_0$ be a block
such that $H \cup B$ \ $Y$-codes $H \cup L_0$.  
Without loss of generality we may assume $(0,-1) \in B$.
Since $F$ intersects the $x$-axis $L$ (for example in $(0,0)$) but contains
no points in $\Z^2 \setminus H$, it follows that an $X$-coloring of $H$ determines
a $Y$-coloring of $H$ and likewise an $X$-coloring of $H \cup L_0$ determines
a $Y$-coloring of $H \cup L_0$.
Let $B_0$ be a block in $L_0$ that contains the block $B$
and blocks on either end of $B$ whose lengths are the diameter of $F$.
Then for any $(i,j) \in B$, we have that $F+(i,j) \in H \cup B_0$.
It follows that an $X$-coloring of
$H \cup B_0$ determines a $Y$-coloring of $H \cup B$.
Since a $Y$-coloring of $H \cup B$ determines a
$Y$-coloring of $H \cup L_0$ which in turn determines
an $X$-coloring of $H \cup L_0$ we have shown that
$H \cup B_0$ \ $X$-codes $H \cup L_0$ and $H$ is closing for $X$.

\end{proof}

\section{Coding corners in closing light cones.}
\subsection{Spacetimes and light cones}
We give a way to extend a one dimensional system to a two dimensional version, with a variant of the definition of a spacetime from~\cite{CFK} (there is also a related notion called the complete history  in Milnor~\cite{M}): 
\begin{defn}
If $X$ is a $\Z^2$-subshift and  $(e_1, e_2)$ is an {\em ordered pair of basis
  vectors} of $\Z^2$ and if the ray spanned by $-e_1$ is
expansive, then $\U = (X, (e_1,e_2))$ is called a {\em spacetime}
and $(e_1,e_2)$ is called its {\em distinguished basis}.
If $\U_1  = (X_1, (e_1,e_2))$ and 
$\U_2  = (X_2, (f_1,f_2))$ are spacetimes, an {\em isomorphism of spacetimes} 
$\Psi\colon \U_1 \to \U_2$ 
is a $\Z^2$-subshift isomorphism $\Psi\colon X_1 \to X_2$ such that
$\Psi \circ T_{e_i} = T_{f_i} \circ \Psi$.
\end{defn}

We require that the negative of the first basis element, $-e_1$, be
expansive (as opposed to $e_1$) for consistency with common usage: 
when $e_1 = (1,0)$ we want
the lower half space of $\R^2$ to code the upper, not the reverse. In addition, 
with this convention  for a polygonal shift with coding polygon $\cP$
we  have that all nonexpansive rays are parallel
(rather than anti-parallel) to edges of $\cP$ with their standard (counterclockwise)
orientation (see, for example, Proposition~\ref{prop: polygon sides}).

Note that this definition of a spacetime $\U$  is more general than that
given in~\cite{CFK}, where it is  required that $ - e_1$ be 
$1$-expansive in the sense that the line $L \cap \Z^2$ containing
$e_1$ codes the half space $\{j e_1 + m e_2 \in \Z^2 \colon m \ge 0\}$. 
This requirement is equivalent to the existence of an endomorphism $\phi$ of a
$\Z$-subshift $\sigma\colon Y \to Y$ with the same alphabet 
as $\U$  such that  $u \in \U$ if and only if 
\begin{enumerate}
\item for each $j \in \Z$ the sequence $\{y_n = u(n, j) \colon n \in \Z\}$ is 
an admissible sequence in $Y$, and
\item if $y \in Y$ satisfies $y_n = u(n, j)$ for all $ n \in \Z$, 
then for $\phi^m(y)_n  = u(n, j+m)$ for all $ n \in \Z$.
\end{enumerate}

When these two conditions are satisfied, we say the spacetime $\U$ is the 
{\em spacetime of the endomorphism $\phi$}. 
We show (see Lemma~\ref{prop: expansive line} below)
that if $\V$ is any spacetime with at least one expansive ray,
then it can be recoded  to be the spacetime of an
endomorphism.

If $\phi\in\End(Y, \sigma)$ and $n\geq 0$, 
following~\cite{CFK} we define $W^+(n, \phi)$ 
to be the smallest integer such that the ray 
$[W^+(n, \phi), \infty)$ is $\phi^n$-coded by $[0,\infty)$ 
and define $W^-(n, \phi)$ to be the largest integer such that the ray 
$(-\infty, W^-(n, \phi)]$ is $\phi^n$-coded by $(-\infty, 0]$.
It is straightforward to check that
\begin{equation}\label{W+ eqn}
W^+(k, \phi \sigma^p) = -pk + W^+(k, \phi) 
\end{equation}
 and 
 $$
W^-(k, \phi \sigma^p) = -pk + W^-(k, \phi),
$$
for all $p \in \Z$  (see~\cite{CFK} for more details, however, note that
the published version of~ \cite{CFK} contains a sign error -- the
negative sign in Equation~\eqref{W+ eqn} is omitted).

\begin{defn}
The {\em future light cone $\C_f(\phi)$ of $\phi \in \End(X)$} is defined to be 
\[
\C_f(\phi) = \{(i,j) \in \Z^2 \colon W^-(j, \phi) \le i \le W^+(j, \phi),\  j \ge 0\}.  
\]
The {\em past light cone $\C_p(\phi)$ of $\phi$} is defined to be  
$\C_p(\phi) = - \C_f(\phi)$.  The {\em full light cone}
$\C(\phi)$ is defined to be $\C_f(\phi) \cup \C_p(\phi)$.
 \end{defn}

We emphasize that $\C_p(\phi)$, the past light cone of $\phi$, is typically {\em not} 
closely related to the light cone of $\phi^{-1}$.

The light cone is naturally stratified into levels: 
define the {\em $n^{th}$ level of $\C(\phi)$} to be the set 
\begin{equation}\label{def I}
\cI(n, \phi):=  \{ i \in \Z \colon (i, n) \in \C(\phi)\}.
\end{equation}

Recall that the edges of a light cone have asymptotic slopes
defined by 
\begin{equation}
\alpha^+ := \lim_{k \to \infty} \frac{W^+(k, \phi)}{k} \label{eqn: alpha}
\end{equation}
and \[
\alpha^- := \lim_{k \to \infty} \frac{W^-(k, \phi)}{k}.
\]
These limits exist by Fekete's Lemma.

The edges of the light cone $\C(\phi)$ are given
by the graphs of the functions $i = W^+(k, \phi)$ 
$i = W^-(k,\phi)$ and have nice asymptotic properties.

\begin{defn} The {\em asymptotic light cone} $A(\phi)$ of $\phi$ 
is defined to be the cone in $\R^2$ 
bounded by the lines
$x = \alpha^+(\phi) y$ and $x = \alpha^-(\phi) y$, meaning that 
\begin{align*}
  A(\phi) = &\{(x,y) \in \R^2 \colon  y \ge 0, \ \alpha^-(\phi) y \le  x \le \alpha^+(\phi)y \}\\
  \cup &\{(x,y) \in \R^2 \colon  y \le 0, \ \alpha^+(\phi) y \le  x \le \alpha^-(\phi)y \}.
\end{align*}
\end{defn}

We view $A(\phi)$ as a subset of $\R^2$ rather
than of $\Z^2$, as we want to consider lines with
irrational slope that may lie in $A(\phi)$ but would intersect
$\C_f(\phi)$ only in $\{0\}$.

The rays 
$t(\alpha^- ,1)$ and $t(-\alpha^+, -1)$ for $t \ge 0$ are nonexpansive rays
(see~\cite[Theorem 4.4]{CFK}), 
where the notation $t(\cdot,1)$ means 
the set of all positive scalar multiples of the vector $(\cdot,1)$.

\begin{defn} We say the asymptotic light cone
$A(\phi)$ has {\em closing edges}
if the rays $t(-\alpha^+ ,-1)$ and $t(\alpha^- ,1)$ are closing (in other words, 
the two rays forming the left edge of $A(\phi)$ are closing).  
\end{defn}

For  $\alpha \in \R$,
define the {\em $\alpha$ quadrants} in $\Z^2$ by
\begin{align*}
Q_1(\alpha) &= \{(i,j) \in \Z^2 \colon j \ge 0, i \ge \alpha j\}\\
Q_2(\alpha)  &= \{(i,j) \in \Z^2 \colon j \ge 0, i \le \alpha j\}\\
Q_3(\alpha)  &= \{(i,j) \in \Z^2 \colon j \le 0, i \le \alpha j\}\\
Q_4(\alpha)  &= \{(i,j) \in \Z^2 \colon j \le 0, i \ge \alpha j\}.
\end{align*}
Even though $Q_i(\alpha)$ is a subset of $\Z^2$ and has no
dependence on any particular spacetime, it  frequently is
the case that we are interested in $Q_i(\alpha)$
as a subset of the domain of colorings in a space time.
This can become confusing when more than one spacetime
is involved.
Hence for clarity we  write 
$Q_i(\alpha, \U)$ to indicate that we are viewing it as a subset of
the domain of the colorings in the spacetime $\U$. 
We refer to a subset of $\Z^2$ as a {\em strip} (respectively,  {\em half
  strip}) if it is the intersection of $\Z^2$ with the set of points in
$\R^2$ between two parallel lines  (respectively,  the intersection of a
strip in $\Z^2$ with a closed half space whose edge is not
parallel to the strip).

\begin{lemma}\label{lem: quad to ray}
Suppose $\U$ is a spacetime, $\alpha$ is rational, 
and the quadrant $Q_4(\alpha)$
$\U$-codes the quadrant $Q_1(\alpha)$.
Then there exists $N>0$ such that the half strip
\[
Q_4(\alpha) \cap  \big(\bigcup_{j = 0}^{N} L_j \big)
\] 
$\U$-codes the quadrant $Q_1(\alpha)$
where $L_j = \{(n, -j) \in \Z^2 \colon n \in \Z\}$ is the horizontal
line in $\Z^2$ through $(0, -j)$.

Moreover there exists a spacetime $\V$ of a $\Z$-subshift endomorphism $\psi$
such that $\V$ is a recoding of $\U$ and such that
the ray $$R := ([0,\infty) \times \{0\}) \cap \Z^2
\subset Q_1(\alpha, \V)$$
$\V$-codes the entire quadrant $Q_1(\alpha, \V)$.
\end{lemma}

\begin{proof}
First assume that $\alpha \ge 0$.
Let $(p,q)$ be the point of $Q_1(\alpha, \V)$ closest to $(0,0)$ such that
$p = \alpha q$ with $q>0$. Hence $\alpha = p/q$ and $p \ge 0$.  
We claim that there exists $N> 0$
such that the finite set of  points  $$T:= \{(r,s) \in Q_1(\alpha, \V) \colon 
0 \le r \le p, \  0 \le s \le q\}$$ 
is  $\U$-coded by the 
half strip  
\[
S(N) := \{(i,j)\colon -N \le j \le 0,\  i \ge \alpha j\}
= Q_4(\alpha) \cap ((-\infty,\infty) \times [-N,0]).
\]
Clearly the points with $s = 0$ are coded since they lie in $S(N)$.
If not all points of $T$ are coded by $S(N)$, 
then there exist $x_n,y_n \in \U$ and $(r_0 , s_0) \in T$  such
that for all $n> 0$ \  $x_n(r_0,s_0) \ne y_n(r_0,s_0)$, but
$x_n$ and $y_n$ agree on the strip $S(n)$.  
By passing to subsequences if necessary, 
we can assume that the sequences $\{x_n\}$ and $\{y_n\}$ converge
to $x_\infty$ and $y_\infty$ respectively with $x_\infty(r_0,s_0) \ne y_\infty(r_0,s_0)$.  
Since these two elements of $\U$ agree on $S(n)$ for all $n$, they agree
on the quadrant $Q_4(\alpha)$. 
This contradicts the hypothesis, proving the claim. 

Since $S(N) + (m,0) \subset S(N)$  for all $m \ge 0$
and the half strip $S(N)$ \  $\U$-codes $T$,  it also $\U$-codes $T + (m,0)$.
Thus it  $\U$-codes  the half strip $S(N) + (p,q)$.
It then follows by induction on $n \ge 0$ that the strip $S(N) + n(p,q)$ 
$\U$-codes $S(N) + (n+1) (p,q)$.
Hence the half strip $S(N)$ \  $\U$-codes the quadrant $Q_1(\alpha)$.
This proves the first assertion of the lemma for $\alpha \ge 0$.

Note that if we define the bi-infinite strip
$$\hat S(N) := (-\infty,\infty) \times [-N,0] = \bigcup_{n \ge 0} (S(N) - (n,0)),$$
then since $S(N)$\ $\U$-codes $Q_1(\alpha)$, we have that
$\hat S(N)$\ $\U$-codes the upper half space $j \ge 0$ and
$\hat S(N)+ (i_0, j_0)$\ $\U$-codes the half space $j  \ge j_0$.

Next we consider $\Psi\colon  (X, \cA) \to (X_F, \cA_F)$, the
canonical recoding of $X$ (see Definition~\ref{def: recoding})
via the finite set $F \subset \Z^2$
which we define to be the triangle 
$$S(N) \cap \{(i,j) \colon i \le 0\} =
\{(i,j) \colon -N \le j \le 0, \  \alpha j \le i \le 0\}.$$
Then $\V := X_F$ is a spacetime and the horizontal axis $j=0$ in $\Z^2$
 $\Psi^{-1}$-codes all horizontal translates of $F$.  But the union of these
horizontal translates is $\hat S(N)$.  Since $\hat S(N)$ \ $\U$-codes
the upper half space $j\ge 0$, it follows that the horizontal axis $j=0$
$\V$-codes the half space $j \ge 0$.  

The ray $R := [0, \infty) \times \{0\}$ (for $\V$)
 $\Psi^{-1}$-codes the half strip $S(N)$ (for $\U$). But $S(N)$
$\U$-codes the quadrant $Q_1(\alpha, \U)$.   
Hence $R$\ $\Psi^{-1}$-codes $Q_1(\alpha ,\U ) \cup S(N)$.
But $Q_1(\alpha, \U ) \cup S(N)$ \  $\Psi$-codes $Q_1(\alpha, \V)$.
Thus $R$\ $\V$-codes $Q_1(\alpha, \V)$.  This completes the
second assertion of the lemma for $\alpha \ge 0$.

The proof is analogous for $\alpha < 0$. 
\end{proof}

\subsection{The role of closing}
Recall that if $H$ is a closing half space with
boundary $L$ and $L_0 = L +z_0$ where $z_0 \in \Z^2$ is chosen
such that $L_0$ is the closest coset of $L$
in the complement of $H$, then there is a finite block $\fB$ in
$L_0 \cap \Z^2$ such that $H \cup \fB$ codes $L_0$. 
We want to show there is a constant $\rho >0$ such that
any $(i,j) \in L_0$, close to $\fB$ is  coded by the set 
$\fB \cup (B_\rho(i,j) \cap H \cap \Z^2)$
where $B_\rho(i,j)\subset \R^2$  is the open ball in $\R^2$  with radius $\rho$.

\begin{lemma}\label{lem: closing-ball}
Let $\ell$ be a rational  closing ray contained in the
one-dimensional subspace $L \subset \R^2$  which bounds the closing
half space $H$. Suppose  $L_0 \ := z+ L, \ z \in \Z^2$, and $\fB \subset L_0$
are as in the definition of closing.  
Suppose further that $(i,j) \in L_0 \cap \Z^2$ and 
$\fB(i,j)$ is a translate of $\fB$ in $L_0 \cap \Z^2$  such that $(i,j) \notin \fB$, but
$(i,j) \in \fB + e_L$, where $e_L$ is a generator of $L \cap \Z^2$.
Then there is a constant $\rho >0$, independent of $i$ and $j$  such that
$\{(i,j)\}$ is coded by the set 
\[
\fB(i,j) \cup (B_{\rho/2}(i,j) \cap H \cap \Z^2), 
\]  where 
$B_{\rho/2}(i,j)$  is the open ball with radius $\rho/2$ centered at $(i,j)$.
\end{lemma}

\begin{proof}
Let $\fB$ be the block  whose existence is guaranteed by the assumption that $\ell$ is closing.
If the result does not hold, then for any
$(i,j) \in L_0$ there exist  sequences $\{x_n\},\ \{y_n\} \in X$
such that $x_n(i,j) \ne  y_n(i,j)$ but $x_n$ and $y_n$ have colorings which agree
on $\fB(i,j) \cup (B_n(i,j) \cap H \cap \Z^2)$.  
Choosing subsequences if necessary we can assume
that there exist $x_\infty, y_\infty\in X$  such that
\[
\lim_{n \to \infty} x_n = x_\infty \text{ and } \lim_{n \to \infty} y_n = y_\infty.
\]
Then the restrictions of $x_\infty$ and $y_\infty$ to  
$\fB(i,j) \cup (H \cap \Z^2)$ are
equal but $x_\infty(i,j)  \ne y_\infty(i,j)$. This contradicts the fact that
$\ell$ is closing and so there exists some value of $\rho$ with the desired property. 
Since such a $\rho$ exists for one $(i,j)$, 
it follows by translating in $L_0$ that the 
same $\rho$ works for any
$(i',j') \in L \cap \Z^2$.  
\end{proof}

While in general  $$W^+(k, \phi) =  \alpha^+ k + \lo(k),$$
it is not in general true that  $W^+(k, \phi) = \lceil \alpha^+ k \rceil$.
However, with appropriate hypotheses we can recode the spacetime of $\phi$
to the spacetime of an endomorphism $\psi$
satisfying $W^+(k, \psi) = \lceil \alpha^+ k \rceil$.
The object of the next three lemmas is to show this holds 
if $\phi$ has a closing light cone.
We begin with some basic facts about $W^+(k, \phi)$ and its
relation to $\alpha^+ k$.  Recall that if $\psi$ is
a recoding of $\phi$ then $\alpha^\pm(\phi) = \alpha^\pm(\psi)$ by
Proposition~5.3 of~\cite{CFK}.

\begin{lemma}\label{ceil ineq}
Suppose $\U$ is the spacetime of an endomorphism $\phi$
and $\V$ is the spacetime of an endomorphism $\psi$
which is a recoding of $\U$ . Then:
\begin{enumerate}
\item\label{it:1} $W^+(k, \phi) \ge \lceil \alpha^+ k \rceil$ for all $k \ge 0$
and $W^-(k, \phi) \le \lfloor \alpha^- k \rfloor$ for all $k \ge 0$.
\item\label{it:2} $W^+(k, \phi) \ge W^+(k, \psi)$ and $W^-(k, \phi) \le W^-(k, \psi)$
for all $k \ge 0$.

\item\label{it:3}
If $W^+(k, \phi) = \lceil \alpha^+ k \rceil$ for all $k \ge 0$,
then $W^+(k, \psi) = \lceil \alpha^+ k \rceil$ for all $k \ge 0$.
Similarly if $W^-(k, \phi) = \lfloor \alpha^- k \rfloor$ for all $k \ge 0$, then
$W^-(k, \psi) = \lfloor \alpha^- k \rfloor, \ k \ge 0$.
\end{enumerate}
\end{lemma}

\begin{proof}
By~\cite[Lemma~4.2]{CFK}, we always
have that $W^+(k, \phi) \ge \alpha^+ k$.  Since $W^+(k, \phi)$ is an integer, 
it follows that $W^+(k, \phi) \ge \lceil \alpha^+ k \rceil$ for all $k \ge 0$.
Similarly $W^-(k, \phi) \le \lfloor \alpha^- k \rfloor$ for $k \ge 0$ 
and so~\eqref{it:1} follows.

To prove~\eqref{it:2}, assume that $F \subset \Z^2$ is finite and
$\Psi \colon  \U \to \V$ is a recoding of $\U$ via $F$.  Let $R(r,s)$ denote the
horizontal $\Z^2$ ray $\{(i,j) \colon i \ge r,\  j =  s\}$.
By the definition of $W^+(n, \phi)$, we have that 
$R(0,0)$ \  $\U$-codes $R( W^+(n, \phi), n)$ for $n \ge 0$. Hence
 $R(i,j)$ \  $\U$-codes $R( i+ W^+(n, \phi), j+n)$ for $n \ge 0$.
 Therefore 
  $\bigcup_{(i,j) \in F} R(i,j)  \text{\ $\U$-codes }$
\[
\bigcup_{(i,j) \in F} R(i+  W^+(n, \phi), j+n)\ 
= \bigcup_{(i,j) \in F+ ( W^+(n, \phi),n)} R(i, j).
\]
But the latter $\Psi$-codes $R(W^+(n, \phi), n)$.  It follows that
$R(0,0)$ \  $\V$-codes $R(( W^+(n, \phi), n)$ for $n \ge 0$.
Thus $W^+(n, \phi) \ge W^+(n, \psi)$.
The fact that $W^-(k, \phi) \le W^-(k, \psi)$ is proved similarly.

To prove~\eqref{it:3}, note that parts~\eqref{it:1} and~\eqref{it:2} imply that 
\[ 
 \lceil \alpha^+ k \rceil \le W^+(k, \psi) \le
W^+(k, \phi) =  \lceil \alpha^+ k \rceil.
\]
The proof for $W^-$ is similar. 
\end{proof}

\begin{lemma}\label{rational edges}
Suppose the asymptotic light cone $A(\phi)$ has closing edges and
that $\alpha^+(\phi) = p/q$ 
and $\alpha^-(\phi) = p'/q$ with $p,p' \ge 0$ and  $q >0$.
Then the spacetime $\U$ of 
$\phi$ can be recoded to
the spacetime $\V$ of an endomorphism of another shift $\psi \in \End(Y)$, for 
which 
\[
W^+(kq, \psi) = \alpha^+ kq
\text{ and } W^-(kq, \psi) = \alpha^- kq
\]
when $k>0$. 
\end{lemma}

\begin{proof}
To prove the equality
$W^+(kq, \psi) = \alpha^+ kq$, it suffices
to consider  the special case $\alpha^+ = 0$.
To see this,  suppose $\alpha^+ = p/q$.
It follows from~\cite[Proposition 3.12]{CFK}, 
or from equations~\eqref{W+ eqn} and \eqref{eqn: alpha}, 
that  $\alpha^+(\psi^m \sigma^{k} ) = -k + m\alpha^+(\psi)$.
Letting $k = p,\ m = q,$ and $\psi' = \psi^q \sigma^{p}$ we have that 
$\alpha^+(\psi') = 0$.
Hence if we show that $W^+(k, \psi') = 0$ for all $k \ge 0$,
then by Equation~\ref{W+ eqn}
\[
0 = W^+(k, \psi^q \sigma^{-p}) = -pk + W^+(k, \psi^q)
= -pk + W^+(kq, \psi),
\]
and so $W^+(kq, \psi) = pk = \alpha^{+}kq$.
Thus it suffices to consider
the special case that there is a recoding $\psi$ of $\phi$
such that $\alpha^+ = \alpha^+(\phi) = \alpha^+(\psi) = 0$.  

Define $\dsp(n) = W^+(n, \phi)  -  \alpha^{+}n$. 
By~\cite[Lemma~4.2]{CFK}, the 
function $\dsp(n)$ is subadditive, nonnegative, and
$\dsp(n) = \lo(n)$.
Since $\alpha^{+} = 0$, it follows that 
$W^+(n, \phi) = \dsp(n)$ for all $n \ge 0$
and we are left with showing that $\dsp(n) = 0$. 

Let $\rho$ be the constant given by Lemma~\ref{lem: closing-ball}.  
Without loss of generality,  we can  assume that $\rho$ is an integer $> 1$.  
Then there exists $C>0$ and arbitrarily large
$n_0$ with the property that if  $r_0 =  \dsp(n_0) +C$, then
\begin{equation}\label{dsp le r}
W^+(k,\phi) = \dsp(k) \le r_0 \quad\text{ for all } 0 \le k \le n_0 
\end{equation}
and
\begin{equation}\label{r_0/n_0}
\frac{r_0}{ n_0} < \frac{1}{2\rho}.
\end{equation}
Namely, to prove~\eqref{dsp le r}, note that if $\dsp(n)$ is bounded for all $n \ge 1$
we can choose $C$ to be an upper bound,
and if $\dsp(n)$ is unbounded we can  choose arbitrarily large $n_0$ such that 
for all $0 \le k \le n_0$, \ $\dsp(k) \le \dsp(n_0)$ and let $C = 0$.  
Then equation~\eqref{r_0/n_0} follows from~\eqref{dsp le r} and the
fact that $\dsp(m) = \lo(m)$.

Observe that if $Q$ is the fourth quadrant
$[0,\infty) \times (-\infty, 0]$, then
\begin{equation}\label{Q eqn}
Q  \quad
\U\text{-codes} \quad  [r_0, \infty) \times [0, n_0],
\end{equation}
where again $r_0 =  \dsp(n_0) +C$.
This holds because $r_0  \ge  \dsp(m) = W^+(m)$
for all $m$ with $0\le m \le n_0$.

But we claim that  also  $([r_0, \infty) \times [0, \infty]) \cup Q$ codes 
$[m, \infty) \times [0, \infty]$ for $0 \le m \le r_0$.  To see this, we 
first  code the vertical line through  $(r_0-1, 0)$ as follows: 
use the one-sided expansiveness of the vertical ray 
$(r_0-1, 0) + t (0, 1),\  t \ge 0$
with a block $$\fB(r_0-1, 0) = \{(m-1, t)\colon  -N \le t \le 0\}.$$
By Lemma~\ref{lem: closing-ball}, we have that $\{(r_0-1, 1)\}$ is coded
by $([r_0, \infty) \times [0, n_0]) \cup Q$ if $n_0 > \rho$.  We can repeat this
using  $\fB(r_0-1, 1) = \{(m-1, t)\colon  -N+1 \le t \le 1\}$
to code $\{(r_0 -1, 2)\}$ and then $\fB(r_0-1, 2) := \{(r_0-1, t)\colon  -N+2 \le t \le 2\}$,
to code $\{(r_0 -1, 3)\}$ etc.  We can continue coding $\{(r_0 -1, k)\}$
so long as $k \le n_0 - \rho$, 
where $\rho> 1$ is the constant from Lemma~\ref{lem: closing-ball}. 

If $H$ is the half space to the right of the line $ (r_0-1 , t),\  t \in \Z$,
then $(r_0 -1, k)$ with $k \le n_0 - \rho$  satisfies
$\fB(r_0 -1, k-1) \cup B_{\rho/2}(r_0 -1, k)  \cap H \cap \Z^2
\subset \fB(r_0 -1, k-1) \cup ([r_0, \infty) \times [0, \infty]) \cup Q$, which codes 
$(r_0 -1, k)$. 
Thus  we have shown that $([r_0, \infty) \times [0, n_0]) \cup Q$ codes 
$[r_0 -1, \infty) \times [0, n_0 - \rho]$.

Since $[r_0, \infty) \times [0, n_0 ] \cup Q$ codes 
$[r_0 -1, \infty) \times [0, n_0 - \rho]$, we can repeat this argument to show that 
 $[r_0 -1, \infty) \times [0, n_0 -\rho] \cup Q$ codes 
$[r_0 -2, \infty) \times [0, n_0 - 2\rho]$, etc.  So,
 as long as $m$ satisfies  $n_0 - m \rho >0$ and $m \le r_0$  we have  that 
\[
[r_0, \infty) \times [0, n_0] \cup Q \text{ codes }
[r_0 - m, \infty) \times [0, n_0 - m \rho ].
\]

But by Equation~\ref{r_0/n_0} above
\[
\frac{r_0}{ n_0} < \frac{1}{2\rho},
\] 
so
$r_0 \rho < n_0/2$. Thus if we take $m = r_0$ then
$n_0 - m \rho = n_0 - r_0 \rho > n_0/2$ so 
\[
[r_0, \infty) \times [0, n_0] \cup Q \text{ codes }
[0, \infty) \times [0, n_0/2 ].
\]

Then by Equation~\ref{Q eqn} we see that
$Q$ codes $[0, \infty) \times [0, n_0/2 ]$.
Since $n_0$ can be arbitrarily large we get that
$Q$ codes the full quadrant  $[0, \infty) \times [0, \infty)$, and in particular the claim follows.  

When $\U$ is the spacetime of the endomorphism $\phi$, it follows from
Lemma~\ref{lem: quad to ray} that there is a spacetime $\V_0$ of
an endomorphism $\psi_0$ and a recoding $\Psi \colon \U \to \V_0$,
with the property that the horizontal ray $[0, \infty) \times \{0\}$
$\V_0$-codes the entire first quadrant of $\Z^2$.  In particular, 
$\V_0$-codes the vertical ray $\{0\} \times [0, \infty)$.  In other
words $W^+(m, \psi_0) = 0$ for all $m \ge 0$ which is
the desired result when $\alpha^+ =0$.
As noted, this suffices to prove the general case  that when 
$\alpha^+(\phi) = p/q$ we have $W^+(kq, \psi_0) = \alpha^+ kq$.

By a similar argument we can recode $\V_0$ to a spacetime $\V$
of the endomorphism $\psi$ with the property that
$W^-(kq, \psi) = \alpha^{-}kq$.  
Since  $\psi$ is a recoding of $\psi_0$,
using part~\eqref{it:3}  of Lemma~\ref{ceil ineq}, it follows that
\[ 
 W^+(kq, \psi) =
W^+(kq, \psi_0) = \alpha^+ kq.
\]
Thus the second recoding did not affect the desired equality
for $W^+$.
\end{proof}

\begin{prop}\label{prop: ceiling}
Suppose $\U$ is the spacetime of the endomorphism $\phi \in \End(Y_0)$ 
whose asymptotic light cone $A(\phi)$ has closing edges and
asymptotic slopes $\alpha^{+} = p/q$ and
$\alpha^{-}= p'/q$. 
Then $\U$  can be recoded to
the spacetime $\V$ of some endomorphism $\psi \in \End(Y_1)$ such that
for all $n \ge 0$,
\[
W^+(n, \psi) =  \lceil \alpha^+ n \rceil \text{  and }
W^-(n, \psi) =  \lfloor \alpha^- n \rfloor.
\]
\end{prop}

\begin{proof}
By Lemma~\ref{rational edges},
after recoding we can assume that there exists $m \ge 0$ such that
\[
W^+(m, \phi) = \alpha^+m \text{ and }W^-(m, \phi)
= \alpha^{-}m.
\]
By~\cite[Lemma~3.10]{CFK}, the function
$W^+(n, \phi)$ is subadditive. 
Hence $W^+(m, \phi) = \alpha^+m$ implies that
$W^+(km, \phi) \le \alpha^+km$ for all $k>0$.
But we always have that $W^+(k, \phi) \ge \alpha^+ k$
(see~\cite[Lemma~4.2]{CFK}) and so $W^+(km, \phi) = \alpha^+km$.

For fixed $i_0, j_0$, define the  ray  
$$R(i_0, j_0) := \{(i,j) \in \Z^2
\colon i \ge i_0, j = j_0\}
$$ to be the positive horizontal ray emanating
from $(i_0, j_0)$. 
Then $R(0, 0)$ codes $R(\alpha^+km, km)$ for all $k \ge 0$.  Translating, 
we obtain  that $R(\alpha^+j,  j)$ codes $R(\alpha^+(j+km), j+ km)$ for all $j \in \Z$.
Hence the half strip $S(m) = \{(i,j)\colon -m \le j \le 0 
\text{ and } i \ge \alpha^+ j\}$ codes 
the quadrant $Q_1(\alpha^+)$.

It follows from the second part of Lemma~\ref{lem: quad to ray}
that $\U$ can be recoded to be the
spacetime $\V_0$ of a $\Z$-subshift endomorphism $\psi_0$
with the property that the horizontal ray $R(0,0)$\   $\V_0$-codes 
the entire quadrant $Q_1(\alpha, \V)$.

Since $(\lceil \alpha^+ j \rceil, j) \in Q_1(\alpha^+)$, 
the ray $[0,\infty)$ \ $\psi_0^j$-codes $[\lceil \alpha^+ j \rceil, \infty)$ and 
so $W^+(j, \psi_0) \le \lceil \alpha^+ j \rceil$.  But
by part~\eqref{it:1} of Lemma~\ref{ceil ineq}, 
we have that $W^+(j, \psi_0) \ge \lceil \alpha^+ j \rceil$
and so $W^+(j, \psi_0) = \lceil \alpha^+ j \rceil$.

We can apply an analogous argument to the spacetime $\V_0$ of $\psi_0$
to obtain a recoding of $\V_0$ to $\V$, the spacetime of an endomorphism
$\psi$ such that $W^-(n, \psi) =  \lfloor \alpha^- n \rfloor$.  This recoding
 still has the property that $W^+(n, \psi) = \lceil \alpha^+ n \rceil$
because part (3) of Lemma~\ref{ceil ineq} asserts
\[ 
 W^+(n, \psi) = 
W^+(n, \psi_0) = \lceil \alpha^+ n \rceil.
\]
Thus the second recoding did not effect the desired equality
for $W^+$.
\end{proof}

Recall that we have defined levels in the light cone of an endomorphism
by $\cI(n, \phi):=  \{ i \in \Z \colon (i, n) \in \C(\phi)\}$.
Our next step is the following lemma about light cones:

\begin{lemma} \label{lem: weakly corner}
Suppose $\phi$ is an endomorphism
of a $\Z$-subshift $(Y,\sigma)$ which has a
closing asymptotic light cone $A(\phi)$ with $\alpha^+ > \alpha^-$. 
Then after recoding, 
there exist integers $m$ and $n_0$ with  $n_0 > m > 0$ such that
$\cI(-n, \phi)$  $\phi^m$-codes $\cI(-n +m, \phi)$  whenever $n \ge n_0$.
\end{lemma}

\begin{proof}
The endomorphism $\phi$ is fixed throughout this proof and so  we simplify
notation by writing $\cI(n)$ for $\cI(n, \phi)$ and 
$W^\pm(n)$ for $W^\pm(n, \phi)$.
Since  the asymptotic slopes $\alpha^+$ and $\alpha^-$ satisfy
$\alpha^+ >\alpha^-$, we have
\[
\lim_{n \to \infty} |\cI(-n)| = \infty, 
\]
where $|\cdot |$ denotes the length of an interval.  

By Lemma~\ref{rational edges}, there exists   $m>0$ 
such that $W^+(jm)/m = j\alpha^+$ and 
$W^-(jm)/m = j\alpha^-$ for all $j >0$.  Indeed $m$ can be chosen to be the
least common multiple of the denominators of $\alpha^+$ and $\alpha^-$.
Also by Proposition~\ref{prop: ceiling} we know that
$|\cI(-n)|$ is monotonically increasing in $n$.

It follows from~\cite[Proposition~3.4]{CFK} that there is a constant 
$C$ such that the interval $[0,C]$ \ $\phi^m$-codes $\{W^+(m)\}$ 
and $[-C, 0]$ \ $\phi^m$-codes $\{W^-(m)\}$. Hence for $t>0$, we have that the 
interval $[0,C+t]$ $\phi^m$-codes $[W^+(m), W^+(m) +t]$ and
the interval $[-C-t, 0]$ $\phi^m$-codes $[W^-(m) -t, W^-(m)]$.
Translating, it follows that for any $t > 0$, we have that $[W^+(-n), W^+(-n) +C +t]$
$\phi^m$-codes $[W^+(-n + m),W^+(-n + m) +t]$.  Therefore
a left-aligned subinterval of $\cI(-n)$ with length $C+t$ 
$\phi^m$-codes a left-aligned subinterval of $\cI(-n+m)$
with length $t$ whenever $t \le |\cI(-n+m)|$ and otherwise 
$\phi^m$-codes all of $\cI(-n+m)$.

Let $t = t(n) := |\cI(-n)| - C$.  Then by
monotonicity of $|\cI(-n)|$, we have that $$t(n) > |\cI(-n+m)| -C.$$ 
Since $|\cI(-n+m)|$ tends to infinity with $n$, there exists
$n_0 > 0$ such that $n \ge n_0$ implies 
\[
t =t(n) \ge |\cI(-n+m)| -C > \frac{|\cI(-n+m)|}{2}.
\]
Thus $\cI(-n)$ codes a left-aligned subinterval of $\cI(-n+m)$
with length $t$ which is greater than half the length of 
$\cI(-n+m)$.  An analogous argument shows that
$\cI(-n)$ codes a right-aligned subinterval of $\cI(-n+m)$
with length  greater than half the length of $\cI(-n+m)$.  
We conclude that  $\cI(-n)$ codes $\cI(-n +m)$ when $n > n_0$.
\end{proof}

\section{Corner coding Sectors}
\label{sec:coding-sectors}

\subsection{Corner coding}
In this section, if $u,v,w \in \R^2$  we write $(u,v)$ or $(u,v,w)$ for the ordered pair
or triple of vectors.  We  say $(u,v,w)$ is  {\em positively cyclically ordered}
if $(u,v)$ and $(v,w)$ are  positively oriented bases of $\R^2.$

\begin{defn}\label{def: sector}
Suppose $\ell_1$ and $\ell_2$ are nonparallel rays in $\R^2$ emanating
from the origin labeled such that for $e_i \ne 0 \in \ell_i$, the basis
$(e_1,  e_2)$ is positively oriented 
and the angle $\gamma$ between $e_1$ and $e_2$ satisfies $0 < \gamma < \pi$.

Define the {\em sector $\cS$ determined by $\ell_1$ and $\ell_2$} to be 
\[
\cS = \Z^2 \cap \big (\ell_1 \cup \ell_2 \cup \{v \in \R^2 \colon (e_1, v, e_2)
\text{ is positively cyclically ordered}\} \big)
\]
for any nonzero $e_1 \in \ell_1, e_2 \in \ell_2$.
The {\em supplementary sector to $\cS$} is defined to be
the sector determined by $\ell_2$ and $-\ell_1$ and is denoted by $\cS_s$.
The sector $\cS$ is {\em rational} if the two rays determining it are rational.
\end{defn}

Note that the sector determined by $\ell_1$ and $\ell_2$ is
the set of points of $\Z^2$ lying either between these rays or on them.

\begin{defn}\label{defn: corner coding}
A rational sector $\cS$ for a $\Z^2$-shift $X$ is {\em corner coding} 
if for any finite set $\F \subset \cS$, the set $\cS \setminus \F$ 
\ $X$-codes all of $\cS$. 
A rational sector $\cS$ for a $\Z^2$-shift $X$ is {\em weakly  corner coding} 
if there is a finite set $\F_0$ such that
for any finite set $\F \subset \cS$ the set $\cS \setminus \F$ 
\ $X$-codes all of $\cS \setminus \F_0$. 
\end{defn}

Equivalently, the rational sector $\cS$  is corner coding if the set $\cS \setminus \{(0,0)\}$ \ $X$-codes
$\{(0,0)\}$ (and hence all of $\cS$). 
This is easily checked by induction on the cardinality of $F$.

We  sometimes make use of 
sectors whose vertex $v$ is not at the origin, for example $\cS = \cS_0 +v$
for some $\cS_0$ a sector as defined in Definition~\ref{def: sector}
with its vertex at the origin.  Extending the definition, we say that 
such an $\cS$ is corner coding if $\cS_0$ is corner coding.  In particular, if
$v$ is a vertex of a polygon $\cP$, then
$\cS(v)$, the {\em sector based at
  the vertex $v$} is defined to be the sector with vertex $v$ and rays emanating
from $v$ containing the edges of $\cP$ which meet at $v$.  
Thus a sector based at a vertex other than the origin is a translate of a sector based at the origin.

\begin{prop}\label{prop: corners=>poly}
Suppose $X$ is a $\Z^2$-subshift and $\cP$ is a convex integer polygon.  
If for each vertex $v\in\cP$  the sector  $\cS(v)$ based at $v$
is corner coding, then for sufficiently large $n> 0$, \ $n\cP$ is a
coding polygon for $X$.  Conversely, if $\cP$ is a coding polygon for
$X$, then for each vertex $v\in\cP$ the sector  $\cS(v)$ based at $v$
is corner coding.
\end{prop}

\begin{proof}
Assume that for each $v\in\cP$, the sector $\cS(v)$ is 
corner coding.  Then there is a finite set $G(v) \subset \cS(v)$ with
$v \notin G(v)$ that codes $v$ (see Lemma~\ref{lem: compactness}).
The sector $n \cS(v)$ with vertex $nv$ has the property that
$nv + G(v)$ codes its vertex $nv$.  It follows that for $n_0$
sufficiently large, any $n \ge n_0 $
satisfies $nv + G(v) \subset n\cP$
and hence $n\cP$ codes $nv$.
Repeating this for each vertex of $\cP$, we obtain $n >0$ such that for
each vertex $w$ of $n \cP, $ the set
$n \cP \setminus \{w\}$ codes $w$.  Hence $n\cP$ is a coding
polygon.

The converse follows immediately from the definition of corner coding.
\end{proof}

\begin{lemma}\label{lem: recoding preserves cc}
If a rational sector $\cS$ is corner coding for $X$ and
$\Psi \colon X \to Y$ is a recoding, then $\cS$ is corner coding for $Y$.
\end{lemma}
\begin{proof}
Suppose $\Psi \colon X \to Y$ is a recoding via the finite set $F$.  By  the equivalent 
formulation of 
Definition~\ref{defn: corner coding}, 
it suffices to show that  the set $\check \cS := \cS \setminus \{(0,0)\}$ 
\ $Y$-codes $\{(0,0)\}$.  But the set $\{(0,0)\}$ \ $\Psi^{-1}$ codes $F$, and so
$\{(i,j)\}$ \ $\Psi^{-1}$ codes $F(i,j) := F + (i,j)$.  Hence
$\check \cS$  \ $\Psi^{-1}$ codes 
 $$\bigcup_{(i,j) \in \check \cS} F(i,j).$$ 
But
\[
\bigcup_{(i,j) \in \check \cS} F(i,j) = \bigcup_ {(r,s) \in F}((r,s) + \check \cS) 
\]
and since $\cS$ is corner coding for $X$, each translate $((r,s) + \check \cS)$ \ $X$-codes
$(r,s)$.  Since this holds for each $(r,s) \in F$, it follows that
$\check \cS$  \ $\Psi^{-1}$-codes $F$.  Thus 
$\check \cS$  \ $Y$-codes $(0,0)$ and hence $\cS$ is corner coding for $Y$.
\end{proof}

\begin{cor}\label{cor: recoding poly}
If $X$ is a polygonal subshift with coding polygon
$\cP$ and $\Psi \colon X \to Y$ is a recoding, then
$Y$ is a polygonal subshift with coding polygon $n\cP$ for some $n>0$.
\end{cor}

\begin{proof}
Since $\cP$ is a coding polygon for $X$, each sector based at a vertex
of the polygon $\cP$ is corner coding for $X$.  By Lemma~\ref{lem: recoding preserves cc}, 
each of these
sectors is corner coding for $Y$.  
Then by Proposition~\ref{prop: corners=>poly}, 
for  sufficiently large $n$, the polygon $n \cP$ is a coding polygon for $Y$.
\end{proof}

However, an example of Salo~\cite{salo-example} shows that 
a system isomorphic to a polygonal system need not itself be polygonal: he constructs a system 
isomorphic to the Ledrappier system with an isomorphism that does not preserve the polygonal property.

\begin{prop}\label{prop: weak=>strong}
If a rational sector $\cS$ is weakly corner coding for $X$, then there is a recoding $Y$ of $X$
for which $\cS$ is corner coding.
\end{prop}

\begin{proof}
By the definition of weakly corner coding, there is a finite set $F_0$ such that
for any finite set $F \subset \cS$, the set $\cS \setminus F$ 
\ $X$-codes all of $\cS \setminus F_0$. Without loss of generality we can assume that
$0 \in F_0$ and $F_0$ is convex.

Choose $\cK$ to be a strip in $\R^2$ with several properties we now describe (see Figure 1).
Assume that $\cK$ crosses both sides of $\cS$ transversely such that
each of its edges intersects the edges of $\cS$ in points of $\Z^2$, and further assume
we choose $\cK$ such that $\cS \setminus \cK$ has
two parts separated by $\cK$: the first $\cB_0$ is finite and the second
$\cB_\infty$ is unbounded.  Assume further that $\cK$ is chosen such that  $\cB_0$
contains the finite set $F_0$.  
Let $\cD = \cK \cap \cS \cap \Z^2$.  Then 
$\cD$ is a finite subset whose convex hull is a trapezoid $\hat \cD$.
Two edges of $\hat\cD$ are antiparallel and lie in
the two edges of $\cK$, and the other two sides of $\hat \cD$ lie in the
two edges of $\cS$.  Note that
$\cD + (i,j) \subset \cS$ for any $(i,j) \in \cS$.
We also assume
that $\cK$ has been chosen  to be sufficiently wide such that 
\[
\cD \cup \cB_\infty = \bigcup_{(i,j) \in \cS} (\cD + (i,j)).
\]
Finally,  let $(m,n)$ be the closest point of $\cD$ to $(0,0)$ and note that
without loss of generality we can assume that $(m,n) \in \ell_1$.  

Note that it suffices to show that the translate $\cS(m,n)$ of $\cS$ is
corner coding.  
\begin{center}
\begin{tikzpicture}
\draw[pattern=north west lines, pattern color=gray] (2.462,-.615)--(3.048,1.143)--(5.333,2)--(4.308,-1.077);
\draw [->]  (0,0)  -- (8,3);
\draw[black,fill=black] (0,0) circle (.3ex);
\node [above] at (-.4,-.1) {$(0,0)$};
\draw [->]  (0,0)  -- (8,-2);
\draw [-]  (2,-2)  -- (4,4);
\draw [-]  (4,-2)  -- (6,4);
\draw [->]  (2.462,-.615)  -- (8,1.462);
\draw[black,fill=black] (2.462,-.615) circle (.3ex);
\node [above] at (8.25,2.7) {$A$};
\node [above] at (8.25,1.25) {$B$};
\node [above] at (8.25,-2.3) {$C$};
\node [above] at (1.7,-1.2) {$(m,n)$};
\node [above] at (3.8,.3) {$\hat \cD$}; 
\node [above] at (4.6,2.8) {$\cK$};  
\node [above] at (1.8,-0.1) {$\cB_0$}; 
\node [above] at (6.5,0.5) {$\cB_\infty$}; 
\draw[pattern=dots, pattern color=gray] (7.7,-1.925)--(4.308,-1.077)--(5.333,2)--(7.7,2.888);
\end{tikzpicture}
\end{center}

\begin{center}
\tiny 
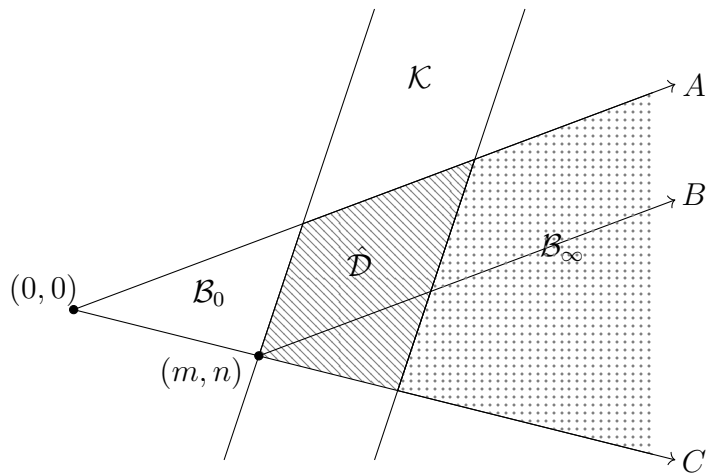
Figure 1: The rays from the origin $(0,0)$ through $A$ and through $C$ are the
edges of the sector $\cS$, and the rays from $(m,n)$ through $B$ and through $C$ are the edges of the sector $\cS(m,n)$.
\end{center}

Set $\cD' = \cD -(m,n)$ and let $\Psi\colon X \to X_{\cD'}$ be the canonical
recoding of $X$ via $\cD'$ to $X_{\cD'}$.  (We use $\cD'$ instead of
$\cD$ because we want $\cD$ to $\Psi$-code $(m,n)$, but using  the canonical
recoding $X_{\cD}$, the set $\cD$ codes $(0,0)$ not $(m,n)$.). 
Thus with the
recoding $\Psi$ via $\cD'$,  we have that $\cD'$ codes $(0,0)$
and so $\cD = \cD' +(m,n)$  \ $\Psi$-codes $(m,n)$.

Thus if $y \in X_\cD$, then the singleton $\{(m,n)\}$ \ $\Psi^{-1}$-codes $\cD$ for the shift $X$.
It follows that the sector $S(m,n)$ \ $\Psi^{-1}$-codes 
$$\bigcup_{(i,j) \in \cS} (\cD+(i,j)) = 
\cD \cup B_\infty$$ for the shift $X$.
Furthermore, $\cD \cup B_\infty$ \  $\Psi$-codes the sector $\cS(m,n)$ for $X_{\cD}$.
Since for any finite set $F$ the set $(\cD \cup B_\infty) \setminus F$ \ $X$-codes 
$\cD \cup B_\infty$, it follows that 
$\cS(m,n) \setminus F$ \  $X_\cD$-codes $\cS(m,n)$ for $X_\cD$.  Thus  
$\cS(m,n)$ is corner coding for $X_\cD$.
Since $\cS(m,n)$ is a translate of $\cS$, it follows that
$\cS$ is corner coding for $X_\cD$.
\end{proof}

\begin{prop}  \label{prop: expansive line}
Suppose $X$ is a $\Z^2$-subshift
and $\C$ is a component of the open set  of expansive rays 
for $X$. Then there exist $u_1, u_2 \in \Z^2$ 
and a $\Z^2$-subshift $Y$ such that: 
\begin{enumerate}
\item $(u_1, u_2)$ is a basis of $\Z^2$ 
\item The rays $\rho_1$ containing $u_1$ and  $\rho_2$ containing  $u_2$ lie in
$\C$.
\item There is a recoding $\Psi\colon X \to Y$.
\item The one-dimensional subspace $L_1$ containing $u_1$ has
the property that $L_1 \cap \Z^2$ codes all of $Y$ and hence it $\Psi^{-1}$-codes all of $X$.
\end{enumerate}
In particular, $Y$ endowed with the basis  $(-u_1, -u_2)$ is 
a spacetime $\U$ of an endomorphism of a $\Z$-subshift. 
\end{prop}

\begin{proof}
Let $\ell_1$ and $\ell_2$ be rays bounding the component $\C$.
Choose a ray $\rho$ in the interior of $\C$ with irrational slope $\lambda$. Let 
$p_1/q_1$ and $p_2/q_2$ be successive convergents for the continued fraction
expansion of $\lambda$ which are chosen such that the subspaces  $L_n$ with
slopes $p_n/q_n$ have slopes sufficiently close to $\lambda$ that the vectors
$u_1 = (p_n, q_n)$ and $u_2 = (p_{n+1}, q_{n+1})$ determine rays $\rho_1, \rho_2$ 
which lie in the interior of  $\C$.  
Since $p_n/q_n$ and $p_{n+1}/q_{n+1}$ 
are successive convergents in the continued fraction expansion of $\lambda$, 
it follows that
$(u_1, u_2)$ is a basis of $\Z^2$  (see Olds~\cite[Section 3.4]{olds}).
Switching the roles of $u_1$ and $u_2$ we can assume it is a positively oriented basis.

Let $L$ be the line containing $\rho_1$.
Since $\rho_1$ is expansive, one of the complementary components 
of $L$ (call this one $H$) codes the other $H'$.  By a change of basis, 
we can assume that $L$ is the horizontal axis and $u_1 = (-1,0)$ and so 
$u_2 = (0,-1)  \in H$. See Figure~2. 

\begin{center}
\begin{tikzpicture}
\draw [->]  (0,0)  -- (-2,0);
\node [left] at (-2,0) {$u_1$};
\draw [->]  (0,0)  -- (0,-2);
\node [below] at (0,-2) {$u_2$};
\draw [dashed]  (0,0)  -- (4,0);
\node [above] at (3.8,0) {$H'$};
\node [below] at (3.8,0) {$H$};
\draw [->]  (4,-2)  -- (-4,2);
\node [above] at (-4,2) {$\ell_1$};
\draw [->]  (-1,3)  -- (1,-3);
\node [below] at (1.4,-2.65) {$\ell_2$};
\node [below] at (2,-1.5) {$\cS_s$};
\draw[black,fill=black] (0,0) circle (.3ex);
\end{tikzpicture}
\end{center}
\begin{center}
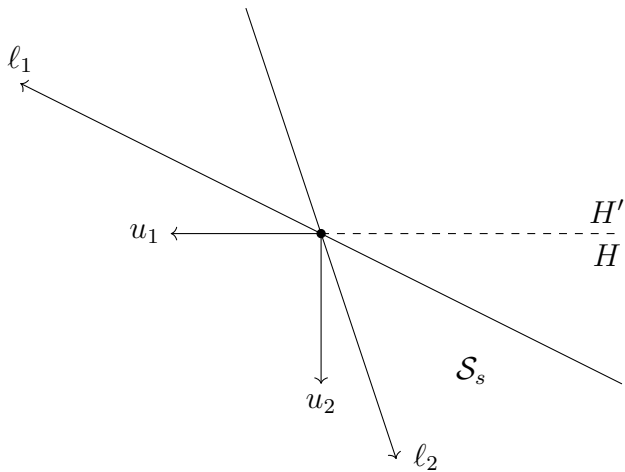
\tiny Figure 2:  The rays $\ell_1$ and $\ell_2$ and the sector $\cS_s$.
\end{center}

Since the negative horizontal axis is an expansive ray,
there exists $r>0$ such that the  strip $S$ consisting
of points of $L \cup H$ with distance at most $r$ from $L$\  $X$-codes the
half space $H'$. 

Let $F$ be the ball in $\Z^2$ of radius $r$ around $ 0$  
and let $Y$ be the shift $X_F$ obtained by the canonical recoding $\Psi\colon X \to X_F$ .
Then $Y$ together with the basis $(-u_1, -u_2)$ is the spacetime of
an endomorphism of the projective subdynamics obtained by restricting $Y$ to $L$.
\end{proof}

\begin{lemma}\label{lem: sector hanging}
Suppose $\cS$ is the sector for the $\Z^2$-subshift $X$
determined by the rays $\ell_1$ and $\ell_2$.
If the supplementary sector $\cS_s$ (bounded by $\ell_2$ and $-\ell_1$) 
is weakly corner coding, then any ray $\ell$ in the interior of the sector $\cS$
is expansive.
\end{lemma}

\begin{proof}
If $L$ is the line containing $\ell$, then $\cS_s \setminus \{(0,0)\}$ lies
in the complementary half space $H$ of  $L$ whose orientation
determines an orientation of $L$ matching that of $\ell$. Let
$H'$ be the other complementary component of $L$.
Translate $\cS_s$ by an element of $\Z^2$ to obtain $\cS_s'$ such that
$B :=\cS_s' \cap H'$ is finite, nonempty, 
and such that some element $b \in B$
is coded by $\cS_s' \cap H$. 
Then by Lemma~\ref{lemma-half}, $H$, and hence
$\ell$, is expansive. See Figure 2.
\end{proof}

\begin{prop}\label{prop: strong coding}
Suppose $\ell_1$ and $\ell_2$ are nonparallel closing rays for
a $\Z^2$-subshift $X_0$ and let 
$\cS_0$ be the sector they determine.
Assume that every ray interior to $\cS_0$ is expansive.  
Let $\cS_s$ be the supplementary sector to $\cS_0$ (the
sector determined by $\ell_2$ and $-\ell_1$).
Then  $X_0$ can be recoded to a $\Z^2$-subshift $X_1$ such that 
the sector $\cS_s$ is corner coding for  $X_1$.
\end{prop}

Note that in this lemma it is not the sector $\cS_0$ which has the corner coding
property, but its supplement $\cS_s$.

\begin{proof}
By hypothesis, the set of all rays in the interior of
$\cS_0$ is a component of the space of expansive rays for $X$.
By Lemma~\ref{prop: expansive line},  after recoding we 
can  assume that $X$ is
the spacetime $(\V, (u_1, u_2))$ of an endomorphism $\psi$ with
$u_1, u_2 \in \cS_0$.
The lines containing the edges of the asymptotic light cone 
$A(\psi)$ of $\psi$  must be the lines containing $\ell_1$ and $\ell_2$,  
since these edges are nonexpansive and there are no other nonexpansive rays 
in $\cS_0$.  Then $\cS_s$ is the lower
half of the asymptotic light cone of this endomorphism. (See Figure 2.)
It follows from  Lemma~\ref{lem: weakly corner}  that
$\cS_s$ is weakly corner coding.
By Proposition~\ref{prop: weak=>strong}, 
there is a recoding such that $\cS_s$ is  corner coding.
\end{proof}

\subsection{Recoding to obtain polygonal shifts}
The primary aim of this section is to prove that if a $\Z^2$-shift X has finitely many nonexpansive rays, all of which are
rational and closing, then X recodes to a polygonal shift.  
Hence, given
a finite set $\E := \{\ell_i\}$ of rational closing rays, which
includes all nonexpansive rays, we want to construct a coding
polygon (for a recoding) whose oriented edges are positively parallel to the
elements of $\E$. Recall (see Section~\ref{sec:parallelisms}) that
an oriented edge $E$ is positively parallel to a ray $\ell$ if a
translate of $E$ lies in $\ell$ with matching orientations, meaning that they
are parallel with matching orientations.

Abstractly, given a set of
rays, a necessary condition for the existence of
a convex polygon with one oriented edge positively parallel to
each ray is that we can find a nonzero vector in each ray such 
that the sum of the vectors is $0$. These vectors are just the
edge vectors of the polygon.  
Thus, given the finite set $\E$ of nonexpansive rays, 
we want  to find a nonzero integer vector $e_i \in \ell_i$ for each $\ell_i \in \E$ 
such that $\sum e_i = 0$ and show these vectors form the edges of a polygon.
We first consider a degenerate case where the polygon
turns out to be a line segment.

\begin{lemma}\label{lem: antiparallel}
Suppose $X$ has two closing and nonexpansive antiparallel rays,
$\ell$ and $-\ell$, which lie in the rational line $L$. If one of the two
components of the complement of $L$ does not intersect any
nonexpansive rays, then $\ell$ and $-\ell$ are the only nonexpansive
rays and the line $L$ determines a periodic direction for $X$.  In
particular, $X$ recodes to a polygonal shift with a degenerate coding
polygon which is a line segment parallel to $L$.
\end{lemma}

\begin{proof}
Without loss of generality, we can assume that that $L$ is vertical, taking $L$ to be the $y$-axis,  and further assume that  the left
half space $H := \{(x,y) \colon x <0\}$ is disjoint from nonexpansive rays.
As we can recode without affecting nonexpansive directions or periodic directions. 
we can do so and further assume that $X$ is the spacetime
of an endomorphism $\phi$ (see Proposition~\ref{prop: expansive line}).
Recall the left upwardly oriented edge of the top half of the 
asymptotic light cone is a nonexpansive ray and there are no other nonexpansive rays
between it and the negative $x$-axis (see ~\cite[Theorem~4.4]{CFK}).
Hence this ray must be either
$\ell$ or $-\ell$; we assume without loss that it is $\ell$.  The same argument shows that
the left downwardly oriented edge of the bottom half of the asymptotic
light cone is  positively parallel to $-\ell$.
It follows that both edges of the asymptotic light cone of $\phi$
must be the line $L$; in other words, the light cone is
degenerate.

Since $\ell$ and $-\ell$ are closing, it follows from Proposition~\ref{prop: ceiling} that
$[0, \infty) \times \{0\}$ codes the first quadrant and hence
$[0, \infty) \times \{1\}$. Likewise 
$(-\infty,0]\times \{0\}$ codes the second quadrant and hence
$(-\infty,0]\times \{1\}$.
Therefore there exists $b>0$ such that for all sufficiently large $c$, the set
$[0, c]\times \{0\}$ codes $[0, c -b] \times \{1\}$ and similarly
$[-c, 0]\times \{0\}$ codes $[- c +b, 0] \times \{1\}$.
Translating the second by $c$, we have that
$[0,c]\times \{0\}$ codes $[b, c] \times \{1\}$.
If $c > 2b$ this implies $[0,c]\times \{0\}$ codes $[0, c] \times \{1\}$, and
so the strip $[0,c]\times [0, \infty)$ is periodic.
It is easy to check that this implies 
$[0,c]\times (-\infty, \infty)$ is periodic.
Hence $L$ is periodic, which implies any non-vertical ray is expansive.
\end{proof}

\begin{lemma}\label{lem: ray set}
Suppose $X$ is an infinite $\Z^2$-subshift with finitely many nonexpansive subspaces, 
all of which are rational and closing.
\begin{enumerate}
\item\label{rayset:one}
If $E \subset \R^2$ is a one-dimensional rational expansive subspace of
$X$, there exist  nonexpansive rays $\rho_1$ and $\rho_2$
and $u_i \ne 0  \in \rho_i,\  i = 1,2$, such that $u_1$ and $u_2$ lie in different components
of  the complement of $E$.  In particular, there are at least $2$ distinct nonexpansive rays.  
\item\label{rayset:two}
If $X$ has only two nonexpansive rays, then they must be antiparallel.
\item\label{rayset:three}
 Suppose  $\{\ell_i\}_{i=1}^n$ 
is the complete set of nonexpansive rays for $X$ and $n \ge 3$. Then  there exist 
nonzero vectors $e_i \in \ell_i \cap \Z^2$ such that $\{e_i\}_{i=0}^n$
(cyclically ordered by angle 
with an axis) are the edges of a convex polygon $\cP$.
\end{enumerate}
\end{lemma}
\begin{proof}
The number of nonexpansive rays is nonzero by~\cite[Theorem 3.7]{BL}, since $X$ is an infinite, compact metric space. 
Part~\eqref{rayset:one} essentially follows from~\cite[Theorem 4.4]{CFK}.
More precisely, if we make a change of basis such that  $E$ is 
horizontal and recode (per Lemma~\ref{prop: expansive line}), then $X$
is the spacetime of an endomorphism $\phi$ and 
the asymptotic light cone $A(\phi)$ of $\phi$ is not empty.
By the same result of~\cite{CFK}, the ray
$\rho_1$ which forms  the left edge of the part of the asymptotic
light cone $A(\phi)$ which
is above the horizontal axis is a nonexpansive ray.  
Similarly the ray
$\rho_2$ which is the left edge of  the part of $A(\phi)$ which
is below the horizontal axis is a nonexpansive ray.
This proves~\eqref{rayset:one}.

To prove~\eqref{rayset:two}, note that if the two rays are not antiparallel
there is an expansive subspace $L$ with the nonzero vectors in
both rays lying on the same side.  This contradicts~\eqref{rayset:one}.

To prove~\eqref{rayset:three}, we first claim that $0$ lies in the interior of the convex hull of
$\{u_i\}_{i=1}^n$ for any choice of $u_i \ne 0 \in \ell_i \cap \Z^2$, and we proceed 
by contradiction.  Recall that $n \ge 3$.
If $0$ does not lie in the interior of the convex hull of
$\{u_i\}_{i=1}^n$  and $u_i \ne 0 \in \ell_i \cap \Z^2$,
then there is a one-dimensional subspace $L$ bounding a  closed, rational half space $H$ 
such that $u_i$ lies in $H$ for all $i$.

If each $u_i$  lies in the interior of $H$, then $L$ is expansive,
a contradiction of~\eqref{rayset:one}.  
If one $u_i$ lies in  $L$ and all the remaining ones lie in
the interior of $H$, then there 
is a subspace $L'$ lying arbitrarily close to $L$ such that all of $\{u_i\}$
lie in the interior of one
component of its complement, again contradicting~\eqref{rayset:one}.

Finally, suppose two of the $u_i$ lie in $L$, and any others lie in the interior
of $H$.  Then the nonexpansive rays containing these two are antiparallel
and by Lemma~\ref{lem: antiparallel}, 
there are no others. Hence we have contradicted the assumption that
there are $n \ge 3$ nonexpansive rays.

This completes the proof of the
claim that $0$ lies in the interior of the convex hull of
$\{u_i\}_{i=1}^n$ for any choice of $u_i \ne 0 \in \ell_i \cap \Z^2$.

We next proceed to the proof of the existence of $\cP$.
Observe that the claim implies that given $\{u_i\}_{i=1}^n$
as above, there are $\{t_i\}_{i=1}^n \subset (0,\infty)$ such that
$$\sum t_i = 1 \text{ and } \sum t_i u_i = 0.$$
Since $u_i \in \Z^2$, all of the $t_i$ can be taken to be  rational.  
Let $e_i = m t_i u_i$ where $m >0$ is chosen such that the vectors $e_i \in \Z^2$.
Then $\sum e_i = 0$ and $e_i \in \ell_i \cap \Z^2$.  We label the $e_i$ such they are cyclically
ordered by the angle they make with the $x$-axis, and form 
a polygonal curve $\cP$ by concatenating translates of the $e_i$ end-to-end
in order.  Since $\sum e_i = 0$, it follows that this defines a closed polygonal curve.
The vertex
where the end of the translate of $e_i$ meets the start of
the translate of $e_{i+1}$ has an exterior angle equal to the
angle between $\ell_i$ and $\ell_{i+1}$.  Since  these exterior angles are all positive
and sum to $2\pi$, it follows  that $\cP$ is a simple closed curve which is convex.  
\end{proof}

\begin{thm}\label{main thm1}
Assume that $X$ is a $\Z^2$-subshift with  a finite nonempty set of nonexpansive 
rays, each of which is rational and closing.  
Then $X$ can be recoded to be a polygonal shift with a polygon $\cP$   having
each oriented edge positively parallel to one of the nonexpansive rays
and each nonexpansive
ray positively parallel to an edge of $\cP$.
\end{thm}

\begin{proof}
By  parts~\eqref{rayset:one} and \eqref{rayset:two} of
Lemma~\ref{lem: ray set}, there are at least
two nonexpansive rays.  If there are exactly two, Lemma~\ref{lem: antiparallel}
implies that $X$ can be recoded to a (periodic) polygonal system with
a degenerate polygon with oriented edges positively parallel to the two expansive rays.

Hence we can assume there are at least three nonexpansive rays.
By part~\eqref{rayset:three} of Lemma~\ref{lem: ray set}, there are vectors
$\{e_i\}$ forming the edges of a polygon $\cT$ with
$e_i \in \ell_i$, where $\ell_i$ denotes the $i^{th}$ nonexpansive ray in $X$ 
and the rays $\{\ell_i\}_{i=1}^n$ are cyclically ordered by the angle made 
with the positive horizontal  axis.   Let $\cS^i$ be the sector
determined by $\ell_i$ and $\ell_{i+1}$.
Note that the angle determined by $\cS^i$ is an exterior angle
of the polygon $\cT$, and there are no nonexpansive rays
in the interior of $\cS^i$.  By Proposition~\ref{prop: strong coding}, we can
recode $X$ such that the supplementary sector $\cS^i_s$ is corner coding.
The sector $\cS^i_s$ is the sector determined by $e_{i+1}$ and $-e_i$, meaning it is
a translate of the $i^{th}$ vertex of $\cT$, and we denote this vertex by by $w_i$.
By repeated recoding, we can guarantee that each $\cS^i_s$ is
corner coding, and it follows from Lemma~\ref{lem: recoding preserves cc} that each
additional recoding does not affect the corner coding properties of previous
corners. Denote the final recoding by $Y$.

By Lemma~\ref{lem: compactness}, there is a finite
set $G_i \subset (\cS^i_s \setminus \{(0,0)\}) $
such that $G_i$ \ $Y$-codes $\{(0,0)\} $ for each $i$.  Setting
$G_i' := G_i + w_i$, we have that $G_i'$ \  $Y$-codes $w_i$.
Choosing $n$ sufficiently large, we can guarantee that the polygon $\cP := n \cT$  contains
$G_i'$ for all $i$.  It follows that for each $i$,  the set
$\cP \setminus \{w_i\}$ \ $Y$-codes $\{w_i\}$, meaning that the polygon
$\cP$ is a coding polygon for $Y$.
\end{proof}

In the spirit  of a converse  to Theorem~\ref{main thm1}, we have: 

\begin{thm}\label{main thm2}
Let $X$ be an infinite polygonal $\Z^2$ subshift.
Assume that $\cP$ is a coding polygon for a recoding $Y$ of $X$ such that $\cP$ has the minimal
number of sides among all coding polygons
for recodings of $X$.
Then each of the oriented edges of $\cP$ determines a ray which
is closing for $X$, and these rays are the only nonexpansive rays for $X$.
\end{thm}

\begin{proof}
By Lemma~\ref{lem: sector hanging}, every subspace not
parallel to an edge of $\cP$ is expansive.   By Proposition~\ref{prop: closing}, oriented edges that are
nonexpansive determine rays which are closing.  It then follows
that every oriented edge determines a  nonexpansive ray
since otherwise by Theorem~\ref{main thm1}
we could produce a recoding with a coding polygon
having fewer sides.
\end{proof}

\begin{cor}\label{cor: min triangle}
Let $X$ be a polygonal $\Z^2$ subshift. 
Then any two minimal recoding polygons for $X$ which are
homothetic differ by a translation.
In particular, if $X$ is triangular, any
two minimal recoding polygons for $X$ differ by a translation.
\end{cor}

\begin{proof}
By Theorem~\ref{main thm2}, any two similar minimal coding polygons $\cP$  and
$\cP'$ have oriented edges which make the same angles with respect
to the axes of $\R^2$.
By translating we can assume that the homothety taking $\cP$ to $\cP'$ fixes
the origin and there are corresponding
edges $e$ and $e'$ which emanate from the origin. Suppose the homothety carrying
$\cP$ to $\cP'$ is multiplication by the rational $r >0$.
If $r>1$, then $\cP$ is a proper subset of $\cP'$ and if $r <1$, then
$\cP'$ is a proper subset of $\cP$.  Hence $r =1$ and $\cP = \cP'$.

If $\cP$  and $\cP'$ are triangles, they must have the same angles and make
the same angles with the axes. It follows that they are homothetic.
\end{proof}

\begin{ex}\label{ex: rect}
We contrast Corollary~\ref{cor: min triangle} with a polygonal system $X$ whose
coding polygon $\cP$ is an $m \times n$  rectangle whose edges are horizontal and vertical.
If we let $\cP'$ be an $(m+1) \times (n+1)$ rectangle containing
$\cP$, then $X$ is polygonal with respect to $\cP'$.
It is easy to check that if $X_F$ is the canonical recoding of $X$
via $F = \cP$, then a $2 \times 2$ square is a coding polygon for $X_F$ and is
the minimal recoding polygon for $X$.
\end{ex}

\section{Directional entropies of polygonal systems}
\label{sec:entropy}

\subsection{Linear polygonal entropy}
We turn to the study of entropy for two dimensional shifts.  
If $X$ is a $\Z^2$-shift with at least one expansive ray, then any finite region 
in the shift is coded by an interval in the expansive direction.  In particular, this 
implies that the two dimensional entropy of any $\Z^2$-shift with at least one expansive ray is zero, 
and so we restrict ourselves to linear entropy.  This leads us to define a generalization of 
directional entropy that depends on a polygon, rather than a line. 
One of the goals of this section is to show that for a polygonal $\Z^2$-subshift with polygon $\cP$, 
there are strong relations between $\cH(X, \cP)$
and the directional entropies.

For a a polygon $\cP$ in $\R^2$ and $r>0$, we denote the $r$-neighborhood  the polygon by $\cP_r$, meaning that 
$$\cP_r = \{ u \in \Z^2 \colon d( u, \cP) < r\}.$$
If $X$ is a $\Z^2$-subshift and $\cS\subset\R^2$, we denote the complexity of $\cS$ in $X$ by $P(X, \cS)$, 
meaning that $P(X, \cS)$ is the number of $X$-colorings of $\cS\cap\Z^2$.  
Milnor~\cite{M} introduced the notion of higher dimensional entropies
(see also~\cite[Section 6]{BL}). We are
interested in the one-dimensional case which we refer to as linear entropy.

\begin{defn}
If $X$ is a $\Z^2$-subshift and $\cP$ is a polygon in $\R^2$, 
define the {\em linear polygonal entropy} of $\cP$ by 
\[
\cH( X, \cP) = \lim_{r \to \infty} \lim_{n \to \infty}\frac{\ln P(X, (n \cP)_r)}{n}.  
\]
\end{defn}
Note that we allow the $\cP$ to be a degenerate polygon, meaning that 
we allow $\cP$ to be a line segment.
If $v \in \R^2$ and $I_v = \{tv \colon t \in [0,1]\}$, then
$\cH( X, I_v)$ is the {\em directional entropy} $h_v(X)$ in the direction $v$ as discussed
by Milnor~\cite{M}.  
Abusing notation slightly, for $v \in \R^2$  we  write 
$\cH( X, v)$ for $\cH( X, I_v)$,
where $I_v$ denotes the
interval $\{t v \colon 0 \le t \le 1\}$ (considered as a 
degenerate polygon in $\R^2$).

If $X = \cA^{\Z^2}$, then $P(X, n\cP)$ is 
exponential in the area of $n\cP$, and thus is an exponential function of something quadratic
that is in $n$.  In particular, this means that in this setting 
$\cH( X, \cP) =\infty$.  On the other hand, this quantity is finite for any $\Z^2$-system with at least one expansive ray (see the remark following Lemma~\ref{H properties}).

We record the following elementary properties of polygonal entropy
(for more details see~\cite[Theorem~6.2]{BL}):
\begin{lemma} \label{H properties}
For a $\Z^2$-subshift $X$, the  polygonal entropy $\cH(X, \cP)$ satisfies  the following properties: 
\begin{enumerate}
\item If $v \in \Z^2$,  the directional entropy $h_v(X)$
corresponding to  $v$ is equal to $\cH( X, v)$.
\item For $v \in \R^2$, $\cH( X, \cP +v ) = \cH( X, \cP)$.
\item For $r >0, \ \cH( X, r\cP) = r\cH( X, \cP)$.  In particular
for $v \in \R^2$,  $\cH( X, rv) = r\cH( X, v)$.
\item\label{H:four} If $\cP_1$ and $\cP_2$ are polygons in $\R^2$ and there are  $v_1, v_2 \in \R^2$ such
that $\cP_1 + v_1 \subset \cP_2 \subset r\cP_1 + v_2$ for some $r \in \Q$, then
\[
\cH(X, \cP_1) \le \cH(X, \cP_2) \le r\cH(X, \cP_1) 
\]
\end{enumerate}
\end{lemma} 

\begin{remark}
If $X$ is a $\Z^2$-subshift with 
at least one expansive ray, then $P(X, n\cP)$ is bounded above by a
linear function of $n$ and so $\cH( X, \cP)$ is finite.  
To see this, note that a long interval $J$
parallel to a one-sided expansive ray  codes a triangle $\cT$ with $J$
on one side.  Taking $J$ to be sufficiently long, then  $\cT$ is large enough
to contain a translate of $\cP$.  It follows from the properties in 
Lemma~\ref{H properties} that 
$$
\cH( X, \cP) \le \cH( X, \cT) = \cH( X, J) < \infty.
$$
\end{remark}

We recall the following definition (see~\cite{MW} for example).
\begin{defn} If $X$ is a $\Z^2$-subshift, the
entropy seminorm for  $X$ on $\R^2$ is defined by 
\[
\|v \|_X = h_v(X).
\]
\end{defn}

In general $\|v \|_X$ defines a seminorm (see~\cite{BL}), 
but  when $X$ is polygonal with respect to 
$\cP$ and no two sides of $\cP$ are antiparallel,  then 
$\|\  \|_X$ is either identically $0$ or a norm.
To prove this, we make use of a small variation of a result of Milnor~\cite{M} (see also Boyle and Lind~\cite[Theorem 6.3, part 4]{BL}):
\begin{lemma}[Milnor~\cite{M}]
\label{lem: entropy norm}
Suppose $X$ is a  $\Z^2$-subshift  with finitely many nonexpansive rays 
and assume that for each ray $\ell  \subset \R^2$, at least one of $\ell$ or $-\ell$
is an expansive ray.
Then the directional entropy $h_v(X)$ is either $0$ for all $v \in \R^2$ or
is nonzero for all  $v \ne 0$.  Thus the entropy seminorm $\|\ \|_X$ is either trivial or a norm.
\end{lemma}

\begin{proof}
A special case of~\cite[part 4, Theorem 6.9]{BL}
implies that the directional entropy function $h_v(X)$ is continuous 
in $v$.  Thus the
set 
$$Z = \{v\in\R^2 \colon \|v\| = 1 \text{ and } h_v(X) = 0\}$$
is a closed subset
of the unit circle $S^1 \subset \R^2$.  
We show that the set $Z$ is also open, and hence is either empty or is all of $S^1$.

Let $v \ne 0$ and let  $J$ be an interval in $\R^2$  that is parallel to $v$ and contains $0$.  
Then for some $n, r > 0$, the set  $(nJ)_r$ codes a rectangle $R$ on one side of 
$nJ$ with two of its edges parallel to $J$. If $h_v(X) = 0$, then 
$\cH(X,J) =0$ implies $\cH(X, R) =0$.  This implies that $\cH( X, I) =0$
for any interval $I$ with endpoints on the ends of $R$ which are perpendicular to $nJ$.
But the unit vectors $v_I$ parallel to such $I$ (with orientation determined
by the orientation $nJ$ inherits from $v$)  form a neighborhood of 
$v$ in $S^1$. Since for each such $v_I$  we have
$h_{v_I} = 0$, it follows that  the set $Z$ is open.
\end{proof}

\begin{cor}\label{cor: poly norm}
Suppose $X$ is a polygonal  $\Z^2$-shift with a coding
polygon having no pairs of antiparallel sides. 
Then the entropy seminorm $\|\ \|_X$ is either trivial or a norm.
\end{cor}

\begin{proof}
By Proposition~\ref{prop: polygon sides}, $X$ 
satisfies the hypotheses of Lemma~\ref{lem: entropy norm} and the statement follows.
\end{proof}

\subsection{Triangular $\Z^2$-systems}

Define the {\em girth  $ \cG( X, \cP)$ of a polygon} $\cP$  in the direction $v$ for some   
nonzero $v \in \Z^2$ to be the maximal length of a line segment that is the intersection of $\cP$ with a line parallel to $v$.
Equivalently, the girth is the smallest distance between two parallel lines which enclose $\cP$ and are orthogonal to $v$.

For triangles, we are able to say more: 
\begin{prop}\label{h_v formula}
Suppose $X$ is a triangular $\Z^2$-subshift  with 
coding polygon $\cT$. 
If $v \in \R^2$, then the directional entropy of $X$ corresponding to a vector $v$ is 
\[
h_v(X) = \frac{ \cH(X, \cT)}{\cG(\cT, v)}\|v\|.
\]
\end{prop}

This means that for a triangular polygonal system $X$,  depending on whether or not  $\cH(X, \cT) = 0$, 
the directional entropy  $h_v(X)$ is either identically $0$ or is nonzero
for all  $v \ne 0.$ 

\begin{proof}
Assume first that the vector $v$ is not parallel to one of the sides of
$\cT$. The girth $\cG(\cT, v)$ is the length of a line segment $J$
parallel to $v$ with one end on a vertex of $\cT$ (which we assume
without loss to be $(0,0)$) and the other end on the side of $\cT$
opposite to this vertex.  The vector from one end of $J$ to the other
is $\cG(\cT, v) \frac {v} {\|v\|}$.  
Hence, 
\[
\cH(X, J) = \frac{\cG(\cT, v)}{\|v\|} h_v(X),
\]  
and so
\[
h_v(X) = \frac{ \cH(X, J)}{\cG(\cT, v)} \|v\|.
\]

We complete the proof by showing that $\cH(X, J) = \cH(X, \cT)$.  The
interval $J$ divides $\cT$ into two smaller triangles $U$ and $W$
which share the common side $J$.  The triangle $U$ shares a vertex $u \ne 0$
with $\cT$ and $W$ shares a vertex $w \ne 0$ with $\cT$.  The sectors
$\cS(u)$ and $\cS(w)$ whose edges are positively parallel to the edges emanating
from $u$ and $w$, respectively, are corner coding sectors for $X$.
Let $L$ be the subspace containing $J$ and choose $r> \diam(\cP)$ so that
all the translates $P(e) := \cP + e$, with $e \in L$, lie  in $L_r$.
Note that for  sufficiently large $n$, 
the union of the translates $P(e)$ which lie in
$nP$  codes all of $nU$. 

Therefore there exists  $n_0 >0$ and $s \in (0,1]$ such that
$(nJ)_r$ $X$-codes $(nU)_{sr}$ for all $n > n_0$.  Similarly, we can 
assume that  $(nJ)_r$ $X$-codes $(nW)_{sr}$ for all $n > n_0$ and
hence $(nJ)_r$ $X$-codes $(n\cT)_{sr}$.  By the definition of
$\cH(X, .)$, we conclude $\cH(X, nJ) \ge \cH(X, n\cT)$ for all $n > n_0.$
Clearly $(nT)_r$ \ $X$-codes $(nJ)_{r}$ and so by the definition of
$\cH(X, .)$ and part~\eqref{H:four} of Lemma~\ref{H properties}, we have
$\cH(X, nJ) \le \cH(X, n\cT)$. Thus $\cH(X, J) = \cH(X, \cT)$.

If $v$ is parallel to one of the sides $J$ of $\cT$, then  that side has length 
$\cG(\cT, v)$.  A similar argument  shows that $\cH(X, J) =  \cH(X, \cT)$, 
and again the result follows.
\end{proof}

\begin{cor}\label{cor: entropy-sphere}
Suppose $X$ is a rational triangular $\Z^2$-subshift 
and suppose $E(\cT)$ is the set of oriented edges of $\cT$. 
If $\cH(X,\cT) \ne 0$, then the unit sphere in the entropy norm $\|\ \|_X$ is
\[
\frac{1}{\cH(X,\cT)} S_X
\]
where $S_X$ is the convex hexagon
whose oriented edges are $\{\pm e \colon e \in E(\cT) \}$.
\end{cor}

\begin{proof}
If $e$ is an oriented edge of $\cT$, then the girth $\cG(\cT, e) = \| e \|$.  Thus by 
Proposition~\ref{h_v formula}, it follows that $h_e(X) = \cH(X, \cT)$. Hence if $e$ is a positively
or negatively oriented edge of $\cT$, then $h_e(X)$ lies on the sphere of radius 
$\cH(X, \cT)$ in the norm $\|\ \|_X$.  Suppose $w_0$ is a vertex of $\cT$ and $e_1, e_2$ are the
edges emanating from $w_0$. If $v$ is a vector from $w_0$ to a point on the opposite side of
$\cT$, then the girth $\cG(\cT,v) = \|v\|$.  So by Proposition~\ref{h_v formula} we have 
$h_v(X) = \cH(X, \cT)$ and hence $v$  lies on the sphere of radius 
$\cH(X, \cT)$ in the norm $\|\ \|_X$. 
\end{proof}

\begin{cor}\label{cor conformal}
Suppose $X, Y$ are triangular $\Z^2$-subshifts  with  nontrivial entropy norms
and assume both are polygonal with respect to the same rational triangle $\cT$.
Then there is a constant $c >0$ such that the entropy norms of $X$ and $Y$ satisfy
$\|v\|_X = c \|v\|_Y $ for all $v \in \R^2$.
\end{cor}

\begin{proof}
Let 
\[
c = \frac{\cH(X, \cT)}{\cH(Y, \cT)}.
\]
The result then follows from Corollary~\ref{cor: entropy-sphere}. 
\end{proof}

Note this implies that if $X, Y$ are polygonal $\Z^2$-subshifts  with respect to 
the same triangle $\cT$, then
the ratio of their directional entropies in the direction $v$ is independent of the choice of $v$.  The 
$\Z^2$-subshifts  $X$ and $Y$ can be different and even have different alphabets, but 
the shape of the unit ball in the entropy norm depends only on the triangle $\cT$ not the shift $X$.

Next we turn to the relationship between entropy norms
$\|\ \|_X$ and $\|\ \|_Y$ when $X$ and $Y$  with
respect to the same polygon $\cP$, but with no restriction on the number of edges in the polygon. 
\begin{prop}\label{prop quasi-conformal}
Suppose $X$ is a polygonal $\Z^2$-subshift with coding polygon $\cP$ and assume that $\cP$ 
has no antiparallel sides. If  $\fF(\cP)$ is the family of all $\Z^2$-subshifts which are polygonal with respect to
$\cP$ and which have nontrivial entropy norms, then 
there is a uniform dilatation constant $D>0$, depending only on $\cP$, which
has the property that 
for all $u, v \in S^1$  we have
\[
\frac{1}{D} \le \frac {h_u(X)}{h_v(X)} \le D
\]
for all $X \in \fF(\cP)$.
\end{prop}
Thus the conclusion means that the entropy norms for elements of $\fF(\cP)$
is a {\em quasi-conformal family} of norms.  

\begin{proof}
Suppose $v$ is a unit vector and let $L_1$ and $L_2$ be the
unique lines parallel to $v$ which intersect $\partial \cP$ and such that
the interior of $\cP$ lies between them.
Since $\cP$ has no antiparallel sides, at least one of these
lines contains intersects $\cP$ only in a vertex $w$.
Without loss, assume that this line is $L_2$.  
Let $T_v$ be the unique triangle such that 
\begin{enumerate}
\item The vertex $w$ of $\cP$ is also a vertex of $T_v$.
\item The two edges of $T_v$ which meet at $w$ contain
the two edges of $\cP$ which meet at $w$.
\item The other two vertices of $T_v$ lie in $L_1$.
\end{enumerate}
 Let $W := W(v)$ be the side of $T_v$ which lies in
$L_1$ and let $|W| = |(W(v))|$ denote the length of $W$.  

Because the vertex $w$ of $\cP$ is corner coding, if we replace $\cP$ by
$n\cP$ for some large $n>0$ (still calling it $\cP$), then for a given $\cP$ and $v$ 
there is $s \in (0,  1]$ such that
$W_r$ \ $X$-codes $(T_v)_{sr}$  and indeed
$(nW)_r$ \ $X$-codes $(nT_v)_{sr}$.
Hence  from the definition of $H$ and  Lemma~\ref{H properties}, 
\begin{align*}
\cH( X, T_v) &= \lim_{r \to \infty} \lim_{n \to \infty}\frac{\ln P(X, (n T_v)_r)}{n}\\
&= \lim_{r \to \infty} \lim_{n \to \infty}\frac{\ln P(X, (n T_v)_{sr})}{n}\\
&\ge  \lim_{r \to \infty} \lim_{n \to \infty}\frac{\ln P(X, (n W)_r)}{n} = H(X, W).
\end{align*}
Since $W \subset T_v$, the reverse inequality holds and so 
$\cH(X, W) = \cH(X, T_v)$.

Choose $K : = K_v >0$ sufficiently large such that a translate of $K \cP$ contains $T_v$.
Then $K$ depends only on $\cP$ and $v$.
Since $\cP \subset T_v$, we have that 
\begin{equation}\label{K(v) eqn}
\cH(X,\cP) \le \cH(X, T_v) \le K(v) \cH(X,\cP).
\end{equation}
Note that $h_v(X) = C_v \cH(X, T_v)$, where $C_v = |(W(v))|$ 
which depends only on $\cP$ and $v$.
Since $h_v(X)$ is Lipschitz in $v$ with a Lipschitz constant independent of
$X$ (see~\cite[part 4, Theorem 6.9]{BL})
there is a neighborhood $N_v$ of $v$ in $S^1$ which is independent
of $X$ and such that for all $u \in N_v$, 
\[
\frac{C_v}{2} \cH(X, T_v) \le h_u(X) \le 2C_v  \cH(X, T_v).
\]
Since $S^1$ is compact, there is a finite subcovering $\{N_{v_i}\}_{i=1}^m$ 
of $\{N_v\}$.  Let $K = \max K(v_i)$.  Then  
Equation~\ref{K(v) eqn} implies that 
\[
\cH(X,\cP) \le \cH(X, T_{v_i}) \le K \cH(X,\cP)
\]
for all $1 \le i \le m$. Setting $C = \max C_{v_i}$ and 
$c = \min C_{v_i}$, then  for each $i$ and $u \in N_{v_i}$ 
\[
\frac{c}{2} \cH(X, \cP) \le \frac{c}{2} \cH(X, T_{v_i}) 
\le h_u(X) \le 2C  \cH(X, T_{v_i}) \le 2C K \cH(X,\cP).
\]
Therefore 
\[
\frac{c}{2} \cH(X, \cP) 
\le h_u(X) \le 2C K \cH(X,\cP)
\] for all $u \in S^1$.

It follows that for all $u_1, u_2$ in $S^1$
\[
\frac{h_{u_1}(X)}{h_{u_2}(X)} 
\le \frac{2CK \cH(X,\cP)}{\frac{c}{2}\cH(X,P)} = \frac{4CK}{c}.
\]
Setting $D := \frac{4CK}{c}$, since $C,K$ and $c$
are independent of $X$, the result follows.
\end{proof}

\section{Further directions}
We have several questions we are unable to answer, and we collect some of these in this section.  
The first is if there is a canonical way to represent a polygonal system:
\begin{question}\label{ques2} If $X$ is an infinite
polygonal system, are minimal recoding polygons for
$X$ unique up to translation? \end{question}
Corollary~\ref{cor: min triangle} shows this holds when $X$ has
a triangular coding polygon and Example~\ref{ex: rect} shows that this holds when $X$ has a
rectangular coding polygon with sides parallel to the axes,   However, even adding an assumption that the polygonal system has no antiparallel sides, we 
can not answer this question.

One of our results has the hypothesis of a coding polygon with no two antiparallel sides,
or equivalently no two nonexpansive rays with opposite directions.  We ask:
\begin{question}\label{ques || sides}
Does the conclusion of Corollary~\ref{cor: poly norm} remain valid without the assumption
of no antiparallel sides?  In other words,  is the entropy seminorm for a polygonal system always
either a norm or trivial?  
\end{question}

There is an example of Hochman~\cite{hochman} which has
exactly two nonexpansive rays, each of which
is the negative of the other and neither of which
is closing. Both this example and its Cartesian product with a polygonal system are not
polygonal (see Proposition~\ref{prop: cartesian}).
It seems likely that there exists an example with
finitely many nonexpansive rays, no two of which
are antiparallel, with at least one of them not closing, but we do not know how to construct such an example.

In Corollary~\ref{cor: entropy-sphere}, we showed that for a triangular system
whose coding triangle $\cT$ has $\cH(X,\cT) \ne 0$, the unit sphere in the entropy 
norm is determined by the triangle $\cT$.  In fact this sphere is
\[
\frac{1}{\cH(X,\cT)} S_X, 
\]
where $S_X$ denotes the convex hexagon
whose oriented edges are $\{\pm e \colon e \in E(\cT) \}$.
\begin{question}\label{ques-entropy-sphere}
Does the analogous result hold for systems that are not necessarily triangular polygonal systems?
More precisely, if $X$ has a minimal recoding polygon $\cP$ with $n$ sides
 and $\cH(X,\cP) \ne 0$,
must the unit sphere in the entropy norm of $X$ be the $2n$-gon
\[
\frac{1}{\cH(X,\cP)} S_X, 
\]
where $S_X$ denotes the convex polygon 
whose oriented edges are $\{\pm e \colon e \in E(\cP) \}$ and $E(\cP)$ denotes the set of oriented edges of $\cP$?  
\end{question}

This question may be easier to answer under the additional hypothesis that $\cP$ has no antiparallel sides.  
A positive answer to this question would imply a positive answer to Question~\ref{ques2} for systems
with $\cH(X,\cP) \ne 0$, meaning that for such systems, the minimal recoding polygon for such an $X$
is unique.

Results of Einsiedler~\cite{E} show that for a large class of
algebraic systems defined over a compact abelian group, including the
Ledrappier example, there are uncountably many invariant subspaces
realizing distinct directional entropies. We ask if this is true in
greater generality.

\begin{question}\label{ques Einsiedler}
Suppose $X$ is a nontrivial
polygonal shift (not necessarily algebraic), and there is a rational direction for which
the directional entropy is positive.
Are there uncountably many closed $\Z^2$-invariant
subspaces of $X$ realizing distinct directional topological entropies in that direction?
Is it possible that all values in an open interval can be realized as
the directional entropies in this direction for closed subsystems of $X$?
\end{question}

While some of our results carry over immediately to dimensions greater than $2$, most of our results depend on the geometry of two dimensions.  More generally, we ask: 
\begin{question}\label{ques hidim}
Taking the obvious generalization of a  polyhedral $\Z^d$-shift for $d\geq 3$ (meaning that 
the coloring of any one vertex of the polyhedron is uniquely determined by the others), 
are there  higher dimensional analogues of our results? 
\end{question}


\begin{thebibliography}{99}
\bibitem{BlM}
{\sc F.~Blanchard \& A.~Maass}. 
Dynamical properties of expansive one-sided cellular automata.
{\em Israel J. Math.} {\bf 99} (1997), 149--174.

\bibitem{BM}
{\sc M.~Boyle  \& A.~Maass}. 
Expansive invertible onesided cellular automata. 
{\em J. Math. Soc. Japan} {\bf  52} (2000), no. 4, 725--740.

\bibitem{BL}
{\sc M.~Boyle \& D.~Lind}. 
Expansive subdynamics. 
{\em Trans. Amer. Math. Soc.} {\bf 349} (1997), no. 1, 55--102


\bibitem{CG}
{\sc C.~F.~Colle \& E.~Garibaldi}. 
An alphabetical approach to Nivat's Conjecture.  {\tt arXiv:1904.04897}

\bibitem{CK}
{\sc V.~Cyr \& B.~Kra}.
Nonexpansive $\Z^2$-subdynamics and Nivat's conjecture. 
{\em Trans. Amer. Math. Soc.} {\bf 367} (2015), no. 9, 6487--6537.

\bibitem{CFK}
{\sc V.~Cyr,  J.~Franks,  \& B.~Kra}.
The spacetime of a shift endomorphism. 
{\em Trans. Amer. Math. Soc.}  {\bf 371} (2019), no. 1, 461--488.  

\bibitem{E}
{\sc M~Einsiedler}.  
Invariant subsets and invariant measures for irreducible actions on zero-dimensional groups.
{\em Bull. London Math. Soc.} {\bf 36} (2004), no. 3, 321--331. 

\bibitem{hochman}
{\sc M.~Hochman}.
nonexpansive directions for $\Z^2$ actions. 
{\em Ergodic Theory Dynam. Systems.} {\bf 31} (2011), no. 1, 91--112. 

\bibitem{KM}
{\sc J.~Kari \& E.~Moutot}.
Nivat's conjecture and pattern complexity in algebraic subshifts. 
{\em Theoret. Comput. Sci.} {\bf 777} (2019), 379--386. 

\bibitem{kitchens}
{\sc B.~Kitchens}. 
Symbolic dynamics. 
One-sided, two-sided and countable state Markov shifts. Universitext. Springer-Verlag, Berlin, 1998.

\bibitem{KS}
{\sc B.~Kitchens \& K.~Schmidt}.
Markov subgroups of $(\Z/2\Z)^{\Z^2}$. 
Symbolic dynamics and its applications (New Haven, CT, 1991), 265--283, 
{\em Contemp. Math.}, {\bf 135}, Amer. Math. Soc., Providence, RI, 1992. 

\bibitem{KS2}
{\sc B.~Kitchens \& K.~Schmidt}.
Mixing sets and relative entropies for higher-dimensional Markov shifts.
{\em Ergodic Theory Dynam. Systems} {\bf 13} (1993), no. 4, 705--735. 

\bibitem{ledrappier}
{\sc F.~Ledrappier}. Un champ markovien puet \^etre d\'entropie nulle et m\'elangeant. {\em C.R. Acad.
Sc. Paris}.  {\bf  287} (1978) 561--563.

\bibitem{LM}
{\sc D.~Lind \& B.~Marcus}.
{\em An introduction to symbolic dynamics and coding.} Cambridge University Press, Cambridge, 1995.

\bibitem{MW}
{\sc R. Miles \& T. Ward}.  
Directional uniformities, periodic points, and entropy. 
{\em Discrete Contin. Dyn. Syst. Ser. B} {\bf 20} (2015), no. 10, 3525--3545. 

\bibitem{M}
{\sc J.~Milnor}. 
On the entropy geometry of cellular automata. 
{\em Complex Systems} {\bf 2} (1988), no. 3, 357--385.


\bibitem{MH}
{\sc M.~Morse \& G.~A.~Hedlund}. Symbolic dynamics II. Sturmian trajectories. {\em Amer. J. Math.}
{\bf 62} (1940) 1--42.

\bibitem{olds}
{\sc C.~D.~Olds}. {\em Continued Fractions.}
MAA, New Mathematical Library {\bf 9}, (1962).


\bibitem{salo-example}
{\sc V.~Salo}.  
A note on directional closing.   {\tt arXiv:1902.02076}


\bibitem{schmidt}
{\sc K.~Schmidt}. Dynamical systems of algebraic origin.  
Birkh\"auser/Springer, Basel, 1995.
\end{thebibliography}
\end{document}